\documentclass[a4paper,notitlepage, 8pt,reqno]{article}
  \usepackage{graphicx}
  \usepackage{mathtext}
 \usepackage[cp1251]{inputenc}
  \usepackage[T2A]{fontenc}
 \usepackage[russian,english]{babel}
 \usepackage[all,arc,poly,2cell,curve,arrow,tips]{xy}
  \usepackage{enumerate, eucal, amsthm,amsmath, amssymb}

  \allowdisplaybreaks[4]

 \theoremstyle{definition}
 \newtheorem{defn}{Definition}

 \theoremstyle{plain}
 \newtheorem{theorem}{Theorem}
 \newtheorem*{thm*}{Theorem}
 \newtheorem{prop}{Proposition}
  \newtheorem*{prop*}{Proposition}
 
  \newtheorem*{cor*}{Следствие}
 \newtheorem{lem}{Lemma}
  \newtheorem*{lem*}{Лемма}
  \newtheorem{conjecture}{Conjecture}

 \theoremstyle{remark}
 
 \newtheorem{remark}[defn]{Remark}
 \newtheorem*{remark*}{Замечание}

  \newcounter{ab}
\title{Introduction to finite   $W$-algebras}

 \author{D. V. Artamonov, Lomonosov Moscow State University\footnote{artamonov.dmitri@gmail.com}}

\begin{document}

 \maketitle

 \begin{abstract}These are notes of lectures given at UN Encuentro 2016 at the Colombia National University. We begin with the definition of infinite $W$-algebras. Then we explain the motivation for the definition if finite $W$-algebras.  Then we present basic facts about the structure and representations of finite $W$-algebras. In these lectures we follow the historical development of the subject.

 \end{abstract}

\tableofcontents

\section{What are $W$-algebras?}

The Virasoro algebra is an algebra of infinitesimal conformal mappings of $\mathbb{C}$ into itself. In the conformal field theory there was discovered a canonical embedding of the Virasoro algebra $Vir$ into  $U'(\widehat{\mathfrak{g}})$ i.e. into  the completed universal enveloping algebra for the affine Lie algebra corresponding to a semisimple Lie algebra $\mathfrak{g}$.
Zamolodchikov in \cite{Z} discovered an example of  an infinite $W$-algebra, which is a certain finitely generated associative subalgebras in  $U'(\widehat{\mathfrak{g}})$  which is an extension of  $U(Vir)$. Later there were found other examples of such extensions. All of them are very difficult, it is actually impossible to write explicitly the generators and defining relations.

But it was discovered in \cite{bai}, \cite{DS}, \cite{Bal} that the classical counterparts of $W$ algebras, which are some Poisson algebras, have a  very simple description: they are reductions of the the canonical Poisson structure on $\widehat{\mathfrak{g}}^*_1$.

Now take $\mathfrak{g}^*$ instead of $\widehat{\mathfrak{g}}^*_1$ and perform a similar reduction. Then obtain a Poisson algebra which is called the classical finite $W$-algebra. This construction was discovered by  Tjin in \cite{t}. In \cite{t1} Tjin and Boer gave a construction of an associative algebra which is a counterpart of this Poisson algebra, this is the finite $W$ algebra.

\section{Infinite $W$-algebras}
Let $\mathfrak{g}$ be a semisimple   Lie algebra and
$\widehat{\mathfrak{g}}$ is the corresponding untwisted affine Lie
algebra.

In the conformal field theory there exists the Sugawara construction which defines elements $L_n$ in
$U'(\widehat{\mathfrak{g}})$,  that give a
realization of the Virasoro algebra $Vir$. Hence we have also an
embedding $U(Vir)\hookrightarrow U'(\widehat{\mathfrak{g}})$.

The construction of the elements $L_n$ is based on the second
Casimir element of the algebra $\mathfrak{g}$. What happens when one
performs this construction using a Casimir element of higher order?

It turn out that in this way one can obtain interesting objects
called {\it $W$-algebras}

\subsection{Affine Lie algebras}

\subsubsection{Invariant bilinear form}

 Fix  a semi-simple Lie algebra $\mathfrak{g}$. One can have in mind as an example the following algebras

 \begin{enumerate}
 \item  $\mathfrak{gl}_n$ - the algebra of all   $n\times n$ matrices over $\mathbb{C}$.

 \item $\mathfrak{sl}_n$ - the subalgebra  $\{x\in \mathfrak{gl}_n,\,\,\,\,tr(x)=0\}$.

  \item  $\mathfrak{o}_n$ - the subalgebra  $\{x\in \mathfrak{gl}_n,\,\,\,\,x^t+x=0\}$.

 \end{enumerate}
   There exists a symmetric non-degenerate invariant 2-form on the algebra  $\mathfrak{g}$, named the Killing form,
it is defined by formula

\begin{equation}
\label{kilfo} B(x,y)=tr(ad_x ad_y),\,\,\, x,y\in \mathfrak{g}
\end{equation}

where $ad_x$ is an linear operator  $\mathfrak{g}\rightarrow
\mathfrak{g}$, defined by the formula

\begin{equation*}
z\mapsto ad_x(z)=[x,z],\,\,\, z\in\mathfrak{g}
\end{equation*}

In the formula \eqref{kilfo} we take a trace of a composition of two
operators of such type. Explicitly this form can be written as
follows.

\begin{enumerate}
\item For  $\mathfrak{gl}_n:\,\,\, B(x,y)=2ntr(xy)-2tr(x)tr(y)$
\item For  $\mathfrak{sl}_n:\,\,\, B(x,y)=2ntr(xy)$
\item  For  $\mathfrak{o}_n:\,\,\, B(x,y)=(n-2)tr(xy)$
\end{enumerate}

\begin{prop}
For a simple Lie algebra this form is positive defined,
non-degenerate. In the case of a simple Lie algebras every other
non-degenerate  invariant bilinear form is proportional to the
Killing form $B(.,.)$.
\end{prop}

Take a base $I_{\alpha}$, $\alpha=1,...,m$ of the Lie algebra  $\mathfrak{g}$. In the examples listed above the natural choice of the base is the following.

\begin{enumerate}
\item For    $\mathfrak{gl}_n$ take  the base  $E_{i,j}$. Then
\begin{equation}
B(E_{i,j},E_{k,l})=2n\delta_{j,k}\delta_{i,l}-2\delta_{i,j}\delta_{k,l}
\end{equation}
\item For    $\mathfrak{sl}_n$ take  the base  $E_{i,j}$ for $i\neq j$ and $E_{i,i}-E_{i+1,i+1}$. Then
\begin{align*}
&B(E_{i,j},E_{k,l})=2n\delta_{j,k}\delta_{i,l},\,\,\, i\neq j,\,\,k\neq l,\\
& B(E_{i,i}-E_{i+i,i+1},E_{k,l})=0,\,\,\, k\neq l,\\
&B(E_{i,i}-E_{i+1,i+1},E_{j,j}-E_{j+1,j+1})=-2\delta_{i+1,j}-2\delta_{j+1,i}.
\end{align*}
\item For    $\mathfrak{o}_n$ take  the base  $F_{i,j}=E_{i,j}-E_{j,i}$. Then
\begin{equation}
B(F_{i,j},F_{k,l})=2n\delta_{j,k}\delta_{i,l}-2n\delta_{i,k}\delta_{j,l}
\end{equation}

\end{enumerate}



\subsubsection{A loop algebra}

Introduce a notation for structure constants of the Lie algebra
$\mathfrak{g}$:

\begin{equation}
[I_{\alpha},I_{\beta}]=f_{\alpha,\beta}^{\gamma}I_{\gamma}.
\end{equation}

In this formula and everywhere below a summation over repeating indices is suggested

\begin{defn}
{\it The loop algebra} $\mathfrak{g}((t))$ associated to
$\mathfrak{g}$, is a Lie algebra generated by elements denoted as
$I_{\alpha}^n$

\begin{align}
\begin{split}
\label{comrelopp}
&[I_{\alpha}^n,I_{\beta}^m]=f_{\alpha,\beta}^{\gamma}I_{\gamma}^{n+m}+n\delta_{n,-m}B(I_{\alpha},I_{\beta})c
\end{split}
\end{align}

\end{defn}

This algebra has a geometric interpretation.

 $$
\mathfrak{g}((t))=\{\text{formal mappings $S^1\rightarrow
\mathfrak{g}$}\},
$$

where $S^1$ is $\{z\in\mathbb{C},\,\,\, |z|=1\}$. If we take
coordinate $z=e^{2\pi i\varphi}$, $\varphi\in [0,1]$ on $S^1$ then a formal mapping is just a formal power series $\sum_{n=-\infty}^{\infty}f_nz^n$.
We can define the following mappings

$$
I_{\alpha}^n=I_{\alpha}z^{-n-1}
$$

Note that since $\mathfrak{g}=<I_{\alpha}>$ an arbitrary mapping
$S^1\rightarrow \mathfrak{g}$ can be represented as follows

$$
f(z)=\sum_{\alpha}F^{\alpha}(z)I_{\alpha},
$$

Let us give an interpretation of the coefficient $F^{\alpha}(z)$. This is a formal power series  $F^{\alpha}(z)=F^{\alpha}_nz^{-n-1}$, whose coefficients  $F^{\alpha}_n$ depend  linearly on $f$. Thus $F^{\alpha}_n$ are values of some linear function  on the loop $f(z)$ that is  $F^{\alpha}_n\in \mathfrak{g}((t))^*$.

These functions have the following description. First of all let us note that since we have a fixed non-degenerate bilinear form $B$ we can identify $\mathfrak{g}$ and $\mathfrak{g}^*$,

If we have a fixed base $I_{\alpha}$ in $\mathfrak{g}$ then we obtain a dual base $I^{\beta}$ (see also section \ref{duab}), such that

$$
B(I^{\beta},I_{\alpha})=\delta_{\alpha,\beta}.
$$

Note that
$$
[I^{\alpha},I^{\beta}]=f_{\gamma}^{\alpha,\beta}I^{\gamma}.
$$

One can take a series $I^{\alpha}(z)=\sum_{n}I^{\alpha}z^{-n-1}$.

\begin{lem}
$F^{\alpha}(z)=I^{\alpha}(z)(f)$.
\end{lem}

\subsubsection{An affine Lie algebra}

We define the affine Lie algebra as a central extension of the loop
algebra.

\begin{defn}
{\it The affine Lie algebra}  $\widehat{\mathfrak{g}}$, associated
with   $\mathfrak{g}$ is a Lie algebra generated by elements
denoted as $I_{\alpha}^n$,  where $n\in\mathbb{Z}$,  and an element
$C$ subject to relations

\begin{align}
\begin{split}
\label{comrel}
&[C,I_{\alpha}^n]=0,\\
&[I_{\alpha}^n,I_{\beta}^m]=f_{\alpha,\beta}^{\gamma}I_{\gamma}^{n+m}+n\delta_{n,-m}B(I_{\alpha},I_{\beta})C
\end{split}
\end{align}
\end{defn}

In other words, the affine Lie algebra is a central extension of
the loop algebra.


\subsection{Currents and Operator Product Expansions}

Let us present a language that comes from the conformal field theory
(see \cite{df}).

\begin{defn}
Take a formal variable $z$ and introduce a formal power series

\begin{equation}
\label{tok} I_{\alpha}(z):=\sum_{n\in\mathbb{Z}}
z^{-n-1}I_{\alpha}^n
\end{equation}

 It is called {\it a current}
\end{defn}

\begin{defn}

Take a formal variable $z$ and introduce a formal power series

\begin{align*}
&Q(z):=\sum_{n\in\mathbb{Z}} z^{-n-k}Q_{k},\\
&Q_{k}=\sum
c^k_{\alpha^n_1...\alpha^n_k}I_{\alpha^n_1}^{n_1}...I_{\alpha^n_k}^{n_k},\,\,\,c^k_{\alpha^n_1...\alpha^n_k}\in
\mathbb{C}
\end{align*}

 It is called {\it a field} of (conformal) dimension $k$
\end{defn}

It is a formal power series whose coefficients belong to
$U'(\widehat{\mathfrak{g}})$. Here $U'(\widehat{\mathfrak{g}})$ is a completion on the usual universal enveloping algebra $U(\widehat{\mathfrak{g}})$. The elements of
$U'(\widehat{\mathfrak{g}})$  are represented as series of type $\sum c_{\alpha_1,...,\alpha_k}^{n_1,...,n_k}I_{\alpha_1}^{n_1}...I_{\alpha_k}^{n_k}$  such that for every fixed $N$ the number of summands including $I_{\alpha}^t$ with $|t|<N$ is finite.

Let us be given  power series  $A(z)=\sum_{n\in\mathbb{Z}}
z^{n}A_{n}$ and $B(w)=\sum_{n\in\mathbb{Z}} w^{n}B_{n}$ whose
coefficients belong to some ring.

\begin{defn}
Take the  product $A(z)B(w)$ and expand it in $z-w$, a coefficient
at $(z-w)^n$ can depend on $w$. Denote a coefficient at $(z-w)^n$ as
$(AB)_n(w)$, thus one has

\begin{equation*}
A(z)B(w)=\sum_{n\in\mathbb{Z}} (z-w)^{n}(AB)_n(w)
\end{equation*}

This expansion is called an {\it operator product expansion}. Shortly we write OPE.
\end{defn}

\begin{defn}
The coefficient $(AB)_{0}(w)$ is called the {\it normal ordered
product}. It is denoted below as   $(AB)(w)$.
\end{defn}





A formula for the singular part of operator product expansion for
two fields carries the same information as the set of commutation
relations between two arbitrary modes of expansions of these fields.

In particular one has the following result.

\begin{theorem}
The commutation relations \eqref{comrel} are equivalent  to the
following relation
\begin{equation}
I_{\alpha}(z)I_{\beta}(w)=\frac{n\delta_{n,-m}B(I_{\alpha},I_{\beta})C}{(z-w)^2}+\frac{f_{\alpha,\beta}^{\gamma}I_{\gamma}(w)}{(z-w)}+...
\end{equation}
\end{theorem}



\subsection{The energy-momentum tensor}

Let us identify the central element $C\in
U'(\widehat{\mathfrak{g}})$ with a complex number $k$, denote the
result of factorization of  $U'(\widehat{\mathfrak{g}})$  by this
relation as  $U'(\widehat{\mathfrak{g}})_k$. We call
$U'(\widehat{\mathfrak{g}})_k$ the completed universal enveloping
algebra {\it of level $k$}.

\subsubsection{The Casimir element of the second order}
\label{duab}
As in the previous section let $I_{\alpha}$ be a base in $\mathfrak{g}$ and let
$I^{\alpha}$ be a dual base. This is another base of $\mathfrak{g}$, such that

\begin{equation}
B(I_{\alpha},I^{\beta})=\delta_{\alpha}^{\beta}.
\end{equation}

For example

\begin{enumerate}
\item In the case $\mathfrak{g}=\mathfrak{gl}_{n}$ if one takes $I_{\alpha}=E_{i,j}$ then $I^{\alpha}=\frac{1}{n}E_{j,i}$
\item In the case $\mathfrak{g}=\mathfrak{o}_{n}$ if one takes $I_{\alpha}=F_{i,j}$ then $I^{\alpha}=\frac{1}{2n}F_{j,i}$
\end{enumerate}

Now define an integer, called {\it the dual Coxeter number}

\begin{equation}
g=f^{\alpha,\beta,\gamma}f_{\alpha,\beta,\gamma}
\end{equation}

Explicitly, one has

\begin{enumerate}
\item For $\mathfrak{gl}_n$  one has $$g=2n$$
\item For $\mathfrak{o}_n$  one has $$g=n$$
\end{enumerate}

\begin{defn}

The element of the universal enveloping algebra

\begin{equation}
\label{c2} T_2=\sum_{\alpha} I_{\alpha} I^{\alpha}\in
U(\mathfrak{g})
\end{equation}

 is central. It is called {\it the Casimir element of the second order}.
\end{defn}
Explicitly, one has

\begin{enumerate}
\item For $\mathfrak{gl}_n$  one has $$T_2=\frac{1}{2n}\sum_{i,j} E_{i,j}E_{j,i}$$
\item For $\mathfrak{o}_n$  one has $$T_2=\frac{1}{n}\sum_{i<j} F_{i,j}F_{j,i}$$
\end{enumerate}

\subsubsection{The energy-momentum tensor}

\begin{defn}
By analogy with the formula \eqref{c2} let us define a field of conformal dimension $2$:

\begin{equation}
\label{tem} T(z)=\frac{1}{k+g}\sum_{\alpha}
(I_{\alpha}I^{\alpha})(z)\in U'(\widehat{\mathfrak{g}})((z))_k
\end{equation}

This  element  is called {\it the energy-momentum tensor}.

\end{defn}
This name
comes from the Wess-Zumino-Witten theory.

\begin{theorem}
One has an OPE

\begin{equation}
\label{opet}
T(z)T(w)=\frac{c/2}{(z-w)^4}+\frac{T(w)}{(z-w)^2}+\frac{\partial
T(w)}{(z-w)}...,\,\,\,\, c=\frac{kdim\mathfrak{g}}{k+g}
\end{equation}
\end{theorem}

Now take a decomposition

\begin{equation}
\label{tl} T(z)=\sum_n L_nz^{-n-2},\,\,\, L_n\in
U'(\widehat{\mathfrak{g}})_k
\end{equation}


Explicitly one has

\begin{equation}
\label{forl} L_n=\sum_{m>0}I_{\alpha}^nI^{\alpha,m-n}+\sum_{m\leq
0}I^{\alpha,m-n}I_{\alpha}^n
\end{equation}

One can find an operator product expansion

\begin{equation}
\label{opet1}
T(z)T(w)=\frac{c/2}{(z-w)^4}+\frac{T(w)}{(z-w)^2}+\frac{\partial
T(w)}{(z-w)}...,\,\,\,\, c=\frac{kdim\mathfrak{g}}{k+g}
\end{equation}

And using the OPE \eqref{opet1}  one can find that the commutation
relations for $L_n$ are the following

\begin{lem}
 \begin{align*}
 &[L_n,L_m]=(m-n)L_{m+n}+\frac{1}{12}(m^3-m)\delta_{m,-n}.
 \end{align*}
\end{lem}

\subsection{The Virasoro algebra and the Sugawara construction}

Consider the differential operators $L_n=z^n\frac{\partial}{\partial
z}$. They satisfy the commutation relations

\begin{equation}
[L_n,L_m]=(m-n)L_{m+n},\,\,\,\, n\in\mathbb{Z}
\end{equation}

This algebra is called the {\it de Witt algebra}. By definition this
is a Lie algebra of vector fields on $\mathbb{C}$.

The Lie algebra of operators $L_n$ has a prominent central extension
called the Virasoro algebra

 \begin{defn}
 The Virasoro algebra is Lie algebra generated by elements $L_n$ and $C$ subject to relations
 \begin{align*}
 & [C,L_n]=0,\\
 &[L_n,L_m]=(m-n)L_{m+n}+\frac{1}{12}(m^3-m)\delta_{m,-n}.
 \end{align*}
 \end{defn}

\begin{equation}
\label{opet1}
T(z)T(w)=\frac{c/2}{(z-w)^4}+\frac{T(w)}{(z-w)^2}+\frac{\partial
T(w)}{(z-w)}...,\,\,\,\, c=\frac{kdim\mathfrak{g}}{k+g}
\end{equation}

One can prove that the OPE \eqref{opet1} is equivalent to the
following commutation relations for $L_n$

\begin{equation}
\label{reall}
[L_n,L_m]=(m-n)L_{m+n}+\frac{1}{12}(m^3-m)\delta_{m,-n}
\end{equation}

 Thus we obtain the following result

 \begin{theorem}[The Sugawara construction]
The elements $L_n$, and defined by the formula \eqref{forl}  and $k$ give an
embedding of the Virasoro algebra into $U'(\widehat{\mathfrak{g}})_k$.

  \end{theorem}

  \begin{remark}
  The Virasoro algebra is infinite dimensional but it
  is finitely generated as a Lie algebra. One can see that it is
  generated by  $k$, $L_{-2}$, $L_{-1}$, $L_0$, $L_1$, $L_2$.
  \end{remark}

\subsection{The Zamolodchikov's $W_3$ algerba} In the paper \cite{Z} Alexander Zamolodchikov
obtained the following natural generalization of the Sugawara
construction. We follow a mathematical paper \cite{B}

 For $\mathfrak{g}=\mathfrak{sl}_n$ he took instead of the a cental element $T_2\in U(\mathfrak{g})$ a third-order
 central element $T_3=d^{\alpha,\beta,\gamma}I_{\alpha}I_{\beta}I_{\gamma}\in U(\mathfrak{g})$ and considered the corresponding field

\begin{equation}
W(z)=d^{\alpha,\beta,\gamma}(I_{\alpha}(I_{\beta}I_{\gamma}))(z).
\end{equation}

Here we must specify the choice of the central element. We choose it
such that the tensor $d^{\alpha,\beta,\gamma}$ is traceless (see
formula \eqref{trl} for the definition of traceless tensors).

Then consider a decomposition

\begin{equation}
W=\sum_{n\mathbb{Z}} W_nz^{-n-3},\,\,\,\, W_n\in
U'(\widehat{\mathfrak{g}})_k
\end{equation}

Explicitly, one has

\begin{align*}
&W_n=d^{\alpha,\beta,\gamma}(\sum_{m>0,k>n}I_{\alpha}^mI_{\beta}^{n-m}I_{\gamma}^{n-m-k}+\sum_{m\leq
0,k>n}I_{\beta}^{n-m}I_{\gamma}^{n-m-k}I_{\alpha}^m+\\&+
\sum_{m>0,k\leq
n}I_{\alpha}^mI_{\gamma}^{n-m-k}I_{\beta}^{n-m}+\sum_{m\leq 0,k\leq
n}I_{\gamma}^{n-m-k}I_{\beta}^{n-m}I_{\alpha}^m)
\end{align*}

Then Zamolodchikov calculated the commutation relations between
these elements and elements $L_n$. For this pupose the OPE was calculated

\begin{align*}
&T(z)W(w)=\frac{3W(w)}{(z-w)^2}+\frac{\partial W(w)}{(z-w)}+...
\end{align*}

From this OPE one gets that
\begin{equation*}
[L_n,W_m]=(2n-m)W_{n+m}.
\end{equation*}

Also he calculated the OPE

\begin{align*}
&W(z)W(w)=\frac{c/3}{(z-w)^6}+\frac{2T(w)}{(z-w)^4}+
\frac{\partial T(w)}{(z-w)^3}+\\
&+\frac{1}{(z-w)^2}(2\beta\Lambda(w)+\frac{3}{10}\partial^2T(w)+R^4(w)))+\\
&+\frac{1}{(z-w)}(\beta\partial\Lambda(w)+
\frac{1}{15}\partial^3T(w)+\frac{1}{2}\partial R^4(w)),
\end{align*}

where

\begin{equation}
\Lambda(w)=(TT)(w)-\frac{3}{10}\partial^2T(w),\,\,\,
\beta=\frac{16}{22+5c}
\end{equation}

and $R^4$ is some field that cannot be expressed through $W_n$ and
$T_m$. But for $k=1$ in $U'(\widehat{\mathfrak{g}})_k$ this field
vanishes!

Then one has the following relations for the elements $W_n$

\begin{align}
\begin{split}
\label{comw}
&[W_m,W_n]=\frac{c}{360}m(m^2-1)(m^2-4)\delta_{m+n,0}+\\
&+(m-n)(\frac{1}{15}(m+n+3)(m+n+2)-\frac{1}{6}(m+2)(n+2))L_{m+n}+\\
&+\beta\Lambda_{m+n},\\
&\Lambda_{m}=\sum_{n}(L_{m-n}L_n)-\frac{3}{10}(m+3)(m+2)L_m
\end{split}
\end{align}

\begin{defn}
The $W_3$ algebra is an associative algebra generated by   the
elements $W_n$, $L_n$, $1$

\begin{align*}
&[1,L_n]=[1,W_m]=0\\
&[L_n,L_m]=(m-n)L_{m+n}+\frac{1}{12}(m^3-m)\delta_{m,-n}\\
&[L_n,W_m]=(2n-m)W_{n+m}\\
&[W_m,W_n]=\frac{c}{360}m(m^2-1)(m^2-4)\delta_{m+n,0}+\\
&+(m-n)(\frac{1}{15}(m+n+3)(m+n+2)-\frac{1}{6}(m+2)(n+2))L_{m+n}+\\
&+\beta\Lambda_{m+n},\\
&\Lambda_{m}=\sum_{n}(L_{m-n}L_n)-\frac{3}{10}(m+3)(m+2)L_m
\end{align*}

\end{defn}

Thus this algebra is realized as an associative subalgebra in
$U'(\widehat{\mathfrak{g}})_1$.

Note that the last relation is non-linear due to the term
$\Lambda_{m}$ thus $W_3$ is not  a universal enveloping of Lie
algebra.

The  algebra $W_3$ contains an associative subalgebra $U(Vir)$.  Usually one say that $W_3$ is an extension of $Vir$.

\begin{remark}
This subalgebra is finitely generated, since it can ge generated by
$C$, $L_{-2}$, $L_{-1}$, $L_0$, $L_1$, $L_2$, $W_0$.
\end{remark}

\begin{remark}
The name "$W$-algebra" comes from the fact that the field $W(z)$ constructed in the present section in  some early papers was denoted just in this paper using the letter "W".
\end{remark}

\subsection{The general definition of a  $W$-algebra}

The general definition of the   $W$-algebra is the following.
Take the completed universal enveloping algebra
$U'(\widehat{\mathfrak{g}})_k$  of level $k$ and a collection of
fields $W^{\alpha}(z)$ whose coefficients belong to
$U(\widehat{\mathfrak{g}})$.
Suppose that

\begin{align*}&[L_n,W_m^{\alpha}]=(m(h^{\alpha}-1)-n)W_{n+m}^{\alpha}\,\,\Leftrightarrow \\&\,\,L(z)W^{\alpha}(w)=\frac{h^{\alpha}W^{\alpha}(w)}{(z-w)^2}+\frac{\frac{d}{dw}W^{\alpha}(w)}{(z-w)}+...,\end{align*}

and the commutator of $W_n$ and $W_m$ can be expressed through
$L_p$, $W_q$ (that is no new elements of $U(\widehat{\mathfrak{g}})$
are needed)

Then we say that  $L_n$, $W_m$ form a    $W$-algebra

\subsection{Examples of $W$-algebras, the algebras $W_N$}

\subsubsection{A definition though higher order Casimirs} The first
papers about  $W$-algebras were written by physicists and they were
devoted to a search of examples of $W$ algebras - see list of
examples in \cite{w}. In particular in this list there are algebras
$W_N$ that direct generalizations of $W_3$.

Let us give a  definition of $W_N$ following \cite{B}.
Take  a Casimir element for the algebra $\mathfrak{sl}_{N}$ of type

\begin{equation}
T_M=d^{\alpha_1,...,\alpha_M}I_{\alpha_1}....I_{\alpha_M},\,\,\,
M=2,...,N,
\end{equation}

where the tensor $d^{\alpha_1,...,\alpha_M}$ is traceless if a  convolution of two arbitrary indices vanishes:

\begin{equation}
\label{trl} d^{\alpha_1,...,\alpha,...,\alpha,....\alpha_M}=0
\end{equation}

One can always chose a central element in such a way.  Consider the
corresponding field

\begin{equation}
T_M(z)=d^{\alpha_1,...,\alpha_M}(I_{\alpha_1}....I_{\alpha_M})(z),\,\,\,
M=2,...,N.
\end{equation}

The normal ordered product is not associative, but since the tensor
$d^{\alpha_1,...,\alpha_M}$ is traceless the placement of brackets
in this product is not essential.

Take it's modes

\begin{equation}
T_M(z)=\sum_{n} z^{-n-M}T_n
\end{equation}

The modes of  $T_M(z)$ $M=2,...,N$  in the case $k=1$ form a
$W$-algebra called the $W_N$ algebra.

However the formulas for commutation relations for the modes of
$T_M(x)$ are too difficult to be written explicitly.

\subsubsection{ A description of $W_N$ using the Miura transformation
}

Nevertheless there are more explicit construction of the
$W_N$-algebra by means of Miura transformation (see \cite{B}).

Take $\epsilon_i$ be the set of set of weight of the vector
representation of $\mathfrak{sl}_{N+1}$, normalized such that
$\epsilon_i\epsilon_j=\delta_{i,j}-\frac{1}{N+1}$. Take simple roots
$\alpha_i=\epsilon_i-\epsilon_{i+1}$.

Take a some abstract fields $\varphi(z)$ such that they satisfy the
relations

\begin{equation}
\frac{d}{dz}\varphi(z)\frac{d}{dw}\varphi(w)=-\frac{1}{(z-w)^2}
\end{equation}

Take a product

\begin{equation}
R=\prod_{j=1}^{N+1}(\alpha_0\frac{d}{dz}-h_j(z)),\,\,\,\,
h_j(z)=i\epsilon_i\partial\varphi(z).
\end{equation}


Consider the decomposition

\begin{equation}
R=-\sum_{i=1}^{N+1}U_i(z)\frac{d^i}{dz^i}
\end{equation}

Then modes of $U_i(z)$ generate the $W_N$ algebra.

\subsubsection{The coset construction}

Let us be given an affine Lie algebra $\widehat{\mathfrak{g}}$ and
it's subalgebra $\mathfrak{h}$. We define an associative subalgebra
in $U'(\widehat{\mathfrak{g}})_k$
\begin{defn}
\begin{align*}
W(\widehat{\mathfrak{g}},\mathfrak{h})=\text{fields with values in
$U'(\widehat{\mathfrak{g}})_k$ that commute with fields  with values in $
U'(\mathfrak{h})_k$}
\end{align*}
\end{defn}

\begin{lem}
When this algebra is finitely generated it is a $W$-algebra.
\end{lem}

For the pair $$\mathfrak{g}\subset \widehat{\mathfrak{g}}$$ the
algebra

$$W(\widehat{\mathfrak{g}},\mathfrak{g})\,\,\,\,\,\text{ for } k=1$$

is the $W_N$ algebra (see \cite{B}).


 Also some recent explicit construction of $W_N$ is given in section \ref{wnm}.

\subsection{ Representations of algebras $W_N$}

The structure of irreducible representations of the Virasoro algebra
is well-understood.

We say that a  representation $V$ of the Virasoro algebra is a
representation with a highest weight  $\lambda$ if there is a vector
$v\in V$, such that

\begin{align*}
&L_nv=0,\,\,\,n>0,\\
&L_0v=\lambda v
\end{align*}

By analogy, one defines a highest weight representation of the
$W_N$-algebra

\begin{align}
&L_nv=0,\,\,\,T^M_nv=0,\,\,\,n>0,\\
&L_0v=\lambda u,\,\,\,T^M_0v=\lambda_Mv,
\end{align}

In the case of  Virasoro algebra  it is known a lot about
representation

\begin{enumerate}
\item Every highest-weight module has a unique proper maximal submodule, the factor by this module is irreducible
\item There exists an explicit criteria which says when the highest-weight module is finite-dimensional.

\item The formula for the character of an irreducible highest-weight
module is know explicitly
\end{enumerate}

For the $W$-algebras the second and third questions were intensively studied, but there were found no good answers \cite{B}.

\subsection{The Casimir algebras}

One can also consider objects close to $W$-algebras - the Casimir
algebras. They are defined as follows.

We take Casimir elements

\begin{equation}
T_M=d^{\alpha_1,...,\alpha_M}I_{\alpha_1}....I_{\alpha_M}
\end{equation}

of the considered Lie algebra and consider the corresponding fields

\begin{equation}
T_M(x)=d^{\alpha_1,...,\alpha_M}(I_{\alpha_1}....I_{\alpha_M})(z),\,\,\,
M=2,...,N
\end{equation}

Then  for an arbitrary $k$ we take a subalgebra in $U'(\widehat{\mathfrak{g}})_k$ generated
by modes of these fields. It is called the Casimir algebra.

\section{Classical infinite $W$ algebras}
In this section we define the classical analogues of $W_N$ algebras. These are Poisson algebras that are "limits"$\,$ of  $W_N$ algebras.  Also one say that these Poisson algebras are "classical analogues" of $W_N$ algebras. Controversially, one says that the $W_N$ algebras are quantization of classical $W_N$ algebras.

It was was discovered that the classical $W_N$ algebras are closely related to integrable systems.
We are going to describe this relations.

\subsection{ Poisson algebras}
\subsubsection{A definition of a Poisson algebra and a Poisson manifold}

A Poisson algebra is an algebra   $A$ equipped with an additional operation called the {\it  Poisson bracket}

 $$\{.,.\}: A\otimes A\rightarrow A,$$

 that satisfies the following properties

 \begin{enumerate}
 \item   $\{f,g\}=-\{g,f\}$,
 \item $\{f+g,h\}=\{f,h\}+\{g,h\}$,  $\{\alpha f,g\}=\alpha\{f,g\},$ where $\alpha\in A$,
 \item $\{f,\{g,h\}\}+\{h,\{f,g\}\}+\{g,\{h,f\}\}=0$,
 \item $\{f,gh\}=\{f,g\}h+g\{f,h\}$
 \end{enumerate}

Shortly, a Poisson algebra is an associative algebra with an additional Lie bracket $\{.,.\}$ which is a derivation with respect to the structure of an associative algebra.

\begin{defn}
A Poisson manifold is a manifold $M$ such that the algebra $C^{\infty}(M)$ has a structure of a Poisson algebra.
\end{defn}
 A nontrivial example is the following:

\begin{align*}
&M=\mathbb{C}^{2n}\text{ with coordinates }p_1,...,p_n,q_1,...,q_n\\
&\{f,g\}=\sum_{i=1}^n\frac{\partial f}{\partial q_i}\frac{\partial g}{\partial p_i}-\frac{\partial f}{\partial p_i}\frac{\partial g}{\partial q_i}
\end{align*}

\subsubsection{Two important examples}

Take the space $\mathfrak{g}^*$,  it has a canonical structure of a Poisson manifold.

Take a base $I_{\alpha}$ in $\mathfrak{g}$. It also is a function on $\mathfrak{g}^*$ which acts on a element $f\in \mathfrak{g}^*$ by the ruler

\begin{equation}
\label{gss}
I_{\alpha}(f):=f(I_{\alpha}).
\end{equation}

 In the case when $\mathfrak{g}$ is finite-dimensional one has

$$
 \mathfrak{g}^{**}=\mathfrak{g},
$$

thus $I_{\alpha}$ is a base in the space  of linear  function on $\mathfrak{g}^*$.

If $[I_{\alpha},I_{\beta}]=f_{\alpha,\beta}^{\gamma}I_{\gamma}$ then

\begin{equation}
\label{kkp}
\{I_{\alpha},I_{\beta}\}=f_{\alpha,\beta}^{\gamma}I_{\gamma}
\end{equation}

But since a Poisson bracket satisfies the Leibnitz ruler this equality defines  the Poisson bracket of two arbitrary polynomials in  $I_{\alpha}$ and then    the Poisson bracket





Actually this construction defines a structure of Poisson manifold on   $\mathfrak{g}^*$ also in the case of arbitrary $\mathfrak{g}$ not necessary finite-dimensional. Indeed one has an embedding $\mathfrak{g}\hookrightarrow \mathfrak{g}^{**}$, thus  a base element  $I_{\alpha}$ can be viewed as an element of  $\mathfrak{g}^{**}$, acting according to the formula \eqref{gss}.

But the image of  $\mathfrak{g}$ is always dense in some sense, hence  the Poisson bracket on   the image of  $\mathfrak{g}$, defined in the formula \eqref{kkp}, can be uniquely continued to the whole $\mathfrak{g}^{**}$.

This construction is called {\it the Kirillov-Kostant Poisson structure  on  $\mathfrak{g}^*$.}

Now take an affine Lie algebra $\widehat{\mathfrak{g}}$, let us describe the explicitly Poisson structure on $\widehat{\mathfrak{g}}^*$.

Take a current $I_{\alpha}(x)$, then it acts on a function on $\widehat{\mathfrak{g}}^*$ by the formula

\begin{equation}
\label{ia}
I_{\alpha}(z)(f)=\sum_n f(I_{\alpha}^n)z^{-n-1}
\end{equation}

\begin{prop}
According to Kirillov-Poisson structure one has

\begin{align}
\label{kps}
\{I_{\alpha}(x),I_{\beta}(y)\}=k\delta'(x-y)B(I_{\alpha},I_{\beta})+f^{\gamma}_{\alpha,\beta}I_{\gamma}(y)\delta(x-y),
\end{align}
\end{prop}























\subsection{A classical limit and a quantization}
\subsubsection{The deformation quantization }
\label{defq}

Let us be given a Poisson  algebra $A$ one can then take the formal power series $A[[h]]$.

\begin{defn}
We say (see \cite{kon}) that a structure $\circ$ of an associative algebra is a quantization of a Poisson algebra structure if for $f,g\in A$
\begin{enumerate}
\item $f\circ d=fg+\sum_{k=1}^{\infty} B_k(f,g)h^k$, where
\item $B_k(f,g)\text{ is a bidifferential operator of order at most } k$.
\item $f\circ d-d\circ f=ih\{f,g\}+O(h^2).$
\end{enumerate}
\end{defn}

\begin{defn}

We say that $A$ is {\it a classical limit} of  $A[[h]]$ and $A[[h]]$ is {\it a quantization} of $A$.

\end{defn}

The following result takes place

\begin{theorem}[\cite{kon}]
An algebra of functions on any finite-dimensional Poisson manifold can be canonically quantized.
\end{theorem}

\subsubsection{Example: the classical Virasoro algebra }

This is a Poisson algebra generated by elements denoted as $L_n$, $n\in\mathbb{Z}$ and an element $1$, subject to relations

\begin{align*}
&i\{L_n,L_m\}=(n-m)L_{m+n}+\frac{1}{12}(n^3-n)\delta_{m+n,0},\\
&\{L_n,1\}=0.
\end{align*}

This algebra is a quantization of the Virasoro algebra in the following sense. Take an $h$-Virasoro algebra, which is a $\mathbb{C}[[h]]$-algebra generated by central element $1$ and $L_n$ subject to relations

\begin{align*}
&[L_n,L_m]=(n-m)h L_{m+n}+\frac{1}{12}(n^3-n)\delta_{m+n,0},\\
&[L_n,1]=0.
\end{align*}

For $h=1$ we obtain the Virasoro algebra

Thus the classical Virasoro algebra is obtained as a classical limit of $h$-Virasoro algebras

\subsubsection{Example: the classical  $W_3$ algebra }


This is a Poisson associative algebra generated by elements denoted as $L_n$, $W_m$ $n,m\in\mathbb{Z}$ and an element $1$, subject to relations

\begin{align*}
&i\{L_n,L_m\}=(n-m)L_{m+n}+\frac{c}{12}(n^3-n)\delta_{m+n,0},\\
&\{L_n,1\}=0,\,\,\,\,\, i\{L_n,W_m\}=(2n-m)W_{m+n},\\
&i\{V_n,V_m\}=\frac{16}{5c}(n-m)\Lambda_{n+m}+(n-m)L_{n+m}\cdot\\&\cdot
(\frac{1}{15}(n+m+2)(n+m
+3)-\\&-\frac{1}{6}(n+2)(m+2)) +\frac{c}{360}n(n^2-1)(n^2-4)\delta_{n+m,0},\\
&\text{Where } \Lambda_{n}=\sum_{m=-\infty}^{\infty} W_{n-m}W_{m}
\end{align*}

In the same manner the classical $W_N$-algebras are defined.

\section{A relation between classical $W_N$ algebras and dual spaces to affine Lie algebras.}

Let us explain two facts

\begin{enumerate}
\item The classical $W_N$-algebras are isomorphic to Poisson algebras of  differential operators (see \cite{Bakas})
\item The  Poisson algebras  of   differential operators can be obtained using Hamiltonian reduction from the  dual space of an affine Lie algebra with Kirillov-Poison  structure (see \cite{DS}, \cite{bai}).
\end{enumerate}

This whole picture was first outlined in \cite{Bal}.

Below we explain the main steps of this construction. Although we can just formulate the conclusion that the  classical $W_N$-algebras are reductions of  Kirillov-Poison  structure we explaint in details these two steps since they give  relation of  classical $W_N$-algebras to integrable systems.

\subsection{Pseudodifferential operators}
\begin{defn}
A pseudodifferential operator is an operator of the from
\begin{equation}
\label{pdl}
L=u_n(z)\partial^n+u_{n-1}(z)\partial^{n-1}+...+u_0(z)+\sum_{k=-1}^{\infty} u_k(z)\partial^{k},\,\,\,\,\,\partial\equiv\frac{d}{dz}
\end{equation}
\end{defn}

The composition of operators $a(z)\partial$ and $b(z)\partial^l $  is defined as follows according to the Leibnitz ruler

\begin{equation}
\label{comos}
a(z)\partial\circ (b(z)\partial^l)= a(z)\partial b(z) \partial^l+a(z)b(z)\partial^{l+1}
\end{equation}

This formula gives  an analogous formula for any positive $k$:

\begin{equation}
\label{comos2}a(z)\partial^k\bullet b(z)\partial^l=\sum_{t=1}^k C_{k}^t a(z) (\partial ^t b(z))\partial^{k-t+l}
\end{equation}
In the case of negative $l$  and positive $k$ these formulas are true by definitions.
To define a composition in the case of negative $k$ by analogy with  the formula \eqref{comos} we need a formula defining the  action of a negative power of $\partial$ onto $b(z)$. This action  is defined as follows

$$
\partial^{-1}b(z)=\sum_{i=0}^{\infty}
(-1)^i(\partial^i b(z))\partial^{-1-i}
$$



Also for \eqref{pdl} denote as $L_{+}$ the  differential operator

\begin{equation*}
L_{+}=u_n(z)\partial^n+u_{n-1}(z)\partial^{n-1}+...+u_0(z)
\end{equation*}

\subsection{An integrable hierarchy}
\subsubsection{A definition of an integrable hierarchy}
Let fix some differential operator $L$ of type

$$L=\partial^n+u_{n-1}(z,t)\partial^{n-1}+...+u_0(z,t)$$

\begin{prop}
There exists a pseudodifferential operator $L^{1/2}$ such that

 $$L^{1/2}L^{1/2}=L$$
\end{prop}

Consider a series equation for $k=1,2,...$

 \begin{equation}
 \label{inh}
 \frac{d}{dt}L=[L_+^{k/2},L]
 \end{equation}

 It can be written explicitly as a compatible infinite system of PDE for $u_i(z,t)$. This system is called {\it an integrable hierachy}.

\subsubsection{Example: the Korteveg de Vries hierarchy}

Take $L=\partial^2+u$ and $k=3$. Then one has explicitly

\begin{align*}
&L^{1/2}=\partial+\frac{u}{2}\partial^{-1}-\frac{u'}{4}\partial^{-2}+o(\partial^{-3})\\
&(L^{3/2})_+=\partial^3+\frac{3}{2}u\partial+\frac{3}{4}u'
\end{align*}

Thus \eqref{inh} is written explicitly as

\begin{equation}
\frac{d}{dt}(\partial^2+u)=\frac{1}{4}u'''+\frac{3}{2}uu'
\end{equation}

Hence we obtain an equation

\begin{equation}
4\frac{d}{dt}u=u'''+6uu'
\end{equation}

This is  the well-known {\it  Korteveg de Vries equation}. If one takes other values of $k$ one obtains other equations of the  KdV hierarchy.

\subsubsection{Example: the Boussinesq equation}

Take  $L=\partial^3+u\partial+v$ and $k=2$. One has

\begin{align*}
&L^{1/3}=\partial+\frac{1}{3}u\partial^{-1}+o(\partial^{-2}),\\
&(L^{2/3})_+=\partial^2+\frac{2}{3}u,
\end{align*}

Thus

\begin{equation}
\frac{d}{dt}(\partial^3+u\partial+v)=(2v'-u'')\partial+v''-\frac{2}{3}u'''-\frac{2}{3}uu'
\end{equation}

Hence we obtain a system of equations

\begin{align*}
&\frac{d}{dt}u=2v'-u''\\
&\frac{d}{dt}v=v''-\frac{2}{3}u'''-\frac{2}{3}uu'
\end{align*}

One can eliminate $v$  and obtain the equation

\begin{equation}
\frac{d^2}{dt^2}u=-\frac{1}{3}u'''-\frac{4}{3}(uu')'
\end{equation}

This is the well-know {\it  Boussinesq equation}.

\subsection{Two Hamiltonian structure associated  with an  hierarchy}

\subsubsection{A Hamiltonian system}

Let us give a definition of a hamiltonian system of equations in the finite-dimensional case.

 Let us be given a Poisson manifold $M$.
Fix a function $H$. Now take local coordinates $x_1,...,x_m$ in some open subspace of $X$ (note that these coordinates are some function on $X$)  and consider a system of equations

\begin{equation}
\frac{d}{dt}x_i=\{x_i,H\}.
\end{equation}

These equation define a flow $x(t)=(x_1(t),...,x_m(t))$ which does not depend on the choice of coordinates. This system is called {\it Hamiltonian} and the triple $(X,\{.,.\}, H)$ is called a {\it Hamiltonian structure}.


 Now let us give an infinite-dimensional generalization. In the finite dimensional case one can write  coordinates $(x_1,...,x_m)$ as a function $u(z)$, $z=1,2,...,m$, such that $u(1)=x_1,...,u(m)=x_m$.
 Thus as an  infinite-dimensional generalization of coordinates   $(x_1,...,x_m)$ we can consider functions $u(z)$, $z\in \mathbb{R}$ or a collection of functions $u_1(z),...,u_n(z)$

Instead of functions  $f(x_1,...,x_n)\in C^{\infty}(X)$  we must take functional $F(u_1,...,u_n)$. Here is an example of a functional

 $$F(u_1,..,u_n)=\int_{\mathbb{R}}f(z,u_i(z),\partial u_i(z),...,\partial^k u_i(z)) dz.$$

For these functionals the term generalized functions is used.
As usual in the theory of generalized functions we can identify a function $u_i(z)$ with a functional $$u(z)\mapsto \int u_i(z)u(z)  dz.$$ The functional of such type are dense in the space of all functionals. Hence to define the Poisson bracket of  two functional it is actually sufficient to define a Poisson bracket of functional $u_i(z)$ and $u_j(w)$. Note that a Poisson bracket of two functional is also a functional.

 Fix a functional $H$.  Instead of a system o ODE we consider a PDE

 \begin{equation}
 \label{hpde}
\frac{d}{dt}u_i(z)=\{u_i,H\}(z).
\end{equation}

 This equation is called {\it Hamiltonian} and the a pair $(\{.,.\}, H)$ is called a {\it Hamiltonian structure}.

\subsubsection{The Hamiltonian structure for an integrable hierarchy}

 The equations of the  hierarchy \eqref{inh} can be represented as Hamiltonian PDE, that is in the form \eqref{hpde}. Let us present a formula for the Poisson bracket and for the functional $H$.

Remind that we are given an operator

$$L=\partial^n+u_{n-1}(z,t)\partial^{n-1}+...+u_0(z,t)$$

Also introduce a functional

 \begin{equation}
 l(v_1,...,v_n)=\int \sum_{i=1}^n u_i(z)v_i(z)
 \end{equation}

Define  the operation $res$ by the formula

$$
res(\sum_{i}X_i\partial^i):=X_{-1}
$$

And the operation $Tr$ as follows

$$
Tr(K)=\int dz\, res(K)
$$

 Put

$$
 U=\partial^{1-n}u_1+\partial^{2-n}u_2+...+\partial^{-1}u_n,
 $$

Then for example

\begin{align*}
&Tr(LV)=Tr(\partial^n+u_{n-1}\partial^{n-1}+...+u_0)(\partial^{1-n}v_2+\partial^{2-n}v_3+...+\partial^{-1}v_n)=\\
&=Tr(...+(u_1v_1+...+u_nv_n)\partial^{-1}+...)=\int dz (u_1v_1+...+u_nv_n)=l(v_1,...,v_n)
\end{align*}




 Since functionals of type $l(v)$ are dense in the space of all functionals to define a Poisson bracket  it is sufficient to define    $\{l_1(v),l_2(w)\}$, where $l_1$ corresponds to the operator $L_1$ and   $l_2$ corresponds to the operator $L_2$.
 \begin{defn}
 We define the Poisson bracket and the Hamiltonian such that
 \begin{equation}
 \{l_1(v),l_2(w)\}=Tr((L_1v)_{+}(L_2w)-(wL_1)(vL_2)_{+}),\,\,\,\, H=Tr(L).
 \end{equation}
\end{defn}


 \subsubsection{The case $L=\partial^2+u$}

Take  $L=\partial^2+u$ and consider equations for $u(z,t)$ that come from the equation \eqref{inh}.


Consider the Poisson bracket $\{u(x),u(y)\}$, repand

$$u(z)=-\frac{6}{c}\sum_{k=-\infty}^{\infty}L_k e^{-ikx}-\frac{1}{4}.$$

Then one has

\begin{equation}
i\{L_n,L_m\}=(n-m)L_{n+m}+\frac{c}{12}(n^3-n)\delta_{n+m,0}.
\end{equation}


Thus we one obtained a classical Virasoro algebra!

\subsubsection{The case $L=\partial^3+u_1\partial+u'_1+u_2$}

Take $L=\partial^3+u\partial+u'+v$, in this case we have the Bousinessq hierarchy.

Explicitely one has then for the second
Poisson structure

\begin{align*}
&\{u(x),u(y)\}=\frac{1}{2}(\partial^3_x+2u\partial_x+u')\delta(x-y),\\
&\{u(x),v(y)\}=\frac{1}{2}(3v\partial_x+2v')\delta(x-y),\\
&\{v(x),v(y)\}=-\frac{1}{6}(\partial^5_x+10u\partial_x+15u'\partial^2_x+9u''\partial_x+\\&+16u^2\partial_x+2u'''+16uu')\delta(x-y).
\end{align*}

Consider decompositions
\begin{align*}
&u(x)=-\frac{12}{c}\sum_{n=-\infty}^{\infty}e^{-inx}L_n+\frac{1}{2},\\
&v(x)=\frac{12}{c}\sqrt{10}\sum_{n=-\infty}^{\infty}e^{-inx}V_n
\end{align*}

Then the modes of these decompositions satisfy

\begin{align*}
&i\{L_n,L_m\}=(n-m)L_{n+m}+\frac{c}{12}(n^3-n)\delta_{n+m,0},\\
&i\{L_n,L_m\}=(2n-m)V_{n+m},\\
&i\{V_n,V_m\}=\frac{16}{5c}(n-m)\Lambda_{n+m}+(n-m)L_{n+m}(\frac{1}{15}(n+m+2)(n+m+3)-\\&-\frac{1}{6}(n+2)(m+2))+\frac{1}{360}(n^2-1)(n^2-4)\delta_{n+m,0},\\
&\Lambda_n=\sum_{m=-\infty}^{+\infty}L_{n-m}L_m
\end{align*}

That is $L_n$, $V_m$ form the classical $W_3$ algebra.

\subsubsection{Higher-order operators}

Take $L=\partial^n+u_{n-2}\partial^{n-2}+...+u_0$.

Then one has

\begin{align*}
&\{u_{n-2}(x),u_{n-2}(y)\}=\frac{1}{2}(\partial^3_x+2u_2\partial_x+u'_2)\delta(x-y)
\end{align*}

Thus if one puts

\begin{align*}
&u_{n-2}(x)=-\frac{12}{c}\sum_{n=-\infty}^{\infty}e^{-inx}L_n+\frac{1}{2},\\
\end{align*}

 then $L_n$ form a classical Virasoro algebra. Thus modes of
 $u_i(x)$ form a extension of the classical Virasoro algebra.

\begin{theorem} The modes of $u_i(x)$ form  a the classical $W_n$ algebra.
\end{theorem}

\section{A reduction of a Poisson structure: an approach based on group actions}

Now we are going to explain the statement that the first Hamiltonian
structure is a reduction of reduction of the Kirillov-Poisson
structure on the dual of an affine Lie algebra.

\subsection{The general construction of a Hamiltonian reduction}

Let us be given a Poisson manifold $M$ such that a Lie group $N$ acts on it. Then for each $h\in N$ we have a vector field $v_{h}(x)$ on   $M$ that describe an  infinitesimal action of $h$. To define it let us write $h=exp(a)$, $a\in Lie N$, then

$$v_{h}(x)=\frac{d}{dt}(exp(ta)x)\mid_{t=0}.$$

Vector fields can be regarded as differentiation of the  algebra of functions  $C^{\infty}(X)$. An operation  \begin{equation}\label{re}\{H,. \}\end{equation} is also  a differentiation of the  algebra of functions $C^{\infty}(M)$.

\begin{defn}
The action  of $N$ on $M$ is called {\it  hamiltonian } if the  vector fields $v_{h}(x)$ are hamiltonian, that is for every   $h\in N$ there exists a function   $\phi_{h}$, such that the vector field can be represented in the form \eqref{re}

$$v_{h}f=-\{\phi_{h},f\}$$
\end{defn}

Now return to the general situation: the group $N$ acts on a Poisson $M$ in a hamiltonian way.
\begin{defn}

The mapping

\begin{align*}
&\mu: M\rightarrow (Lie H)^* \text{ such that } <\mu(x),a>=\phi_h(x)\text{  where } h=ext(a).
\end{align*}
 is called {\it a momentum map}.

\end{defn}

 Take an orbit $O$ of coadjoint action of $G$ on $(Lie N)^*$ and a submanifold $\mu^{-1}(O)\subset M$. The key fact is the following.

 \begin{prop} The action of $N$ on $M$ preserves  $\mu^{-1}(O)$ \end{prop}

\begin{equation}
\label{hr}
M_O:=\mu^{-1}(O)/N
\end{equation}

This manifold is called a {\it Poisson reduction of $M$ with respect to action of $N$}.

The functions on $M_O$ can be viewed as functions on $\mu^{-1}(O)$ invariant under the action of $N$.

\begin{lem}
The functions on  $\mu^{-1}(O)$ invariant under the action of $N$ form a Poisson subalgebra in $\mathbb{C}^{\infty}(M)$.
\end{lem}
 This lemma is very easy. We need to prove that in $f$ and $g$ are invariant function then $\{f,g\}$ is also an invariant function. But since the action of $N$ is hamiltonian the fact that $f$ and $g$ are invariant under the action of $N$ is written as follows

 $$
 \{\phi_h,f\}=0,\,\,\, \{\phi_h,g\}=0\,\,\,\forall h\in LieN
 $$

Put

$$
\{\phi_h,\{f,g\}\}=-\{g,\{\phi_h,f\}\}-\{f,\{g,\phi_h\}\}=0
$$

Thus  $\{f,g\}$ is also an invariant function.

As a corollary we obtain a natural Poisson bracket on the algebra of  $N$-invariant funcitons.

Thus we obtain the following result.

\begin{theorem}
 $M_O$ is a Poisson manifold
\end{theorem}

\subsection{The main example}

The main example of a hamiltonian action for us will be the following. Take our affine Lie algebra $\widehat{\mathfrak{g}}$, consider it's dual space  $\widehat{\mathfrak{g}}^*$.

As a manifold we take the hyperplane $$\widehat{\mathfrak{g}}^*_1:=\{f\in \widehat{\mathfrak{g}}^*:\,\,\, f(C)=1\}.$$

In the considered example the momentum mapping is just an embedding

\begin{equation}
\mu: \widehat{\mathfrak{g}}^*_1\hookrightarrow \widehat{\mathfrak{g}}^*
\end{equation}

Thus  \begin{equation}M_O=O\end{equation} is set of orbits of coadjoint action that are contained in  $\widehat{\mathfrak{g}}^*_1$.

Take $\mathfrak{g}=\mathfrak{gl}_{n}$. Take the  loop group $N(z)$ of  upper triangular matrices with units on the diagonal. It acts in an adjoint way on $\widehat{\mathfrak{g}}^*_1$.

That is the elements of $N(z)$ are matrices $$\begin{pmatrix}  1 & n_{1,1}(z) &...& n_{1,n}(z)\\
0& 1 &...& n_{2,n}(z)\\...\\ 0&0&...&0 \end{pmatrix},$$

where $n_{i,j}(z)=\sum_{k\in \mathbb{Z}}n_{i,j}^kz^k$. Take  as  $O$ the $N(z)$-orbit of the element

$$
e=\begin{pmatrix}
0&1&0&...&0&0\\
0&0&1&...&0&0\\
...\\
0&0&0&...&0&1\\
0&0&0&...&0&0\\
\end{pmatrix}
$$

\begin{theorem}
The first Hamiltonian structure is a reduction of  $\widehat{\mathfrak{gl_n}}^*_1$ by the action of the group  $N(z)$.
\end{theorem}

\section{Why affine Lie algebras are related to differential operators?}

The key fact is the following  well-know theorem

\begin{theorem}
The space $\widehat{\mathfrak{g}}^*_1$ is naturally isomorphic to space of connections i.e. differential operators

\begin{equation}\label{do}\frac{d}{dz}+A(z),\end{equation} where $A(z)$ is a $\mathfrak{g}$-valued function, thus this is an element of  $\mathfrak{g}((z))\simeq \widehat{\mathfrak{g}}^*_1 $

The coadjoint action of $GL(z)$ corresponds to the  gauge group action of this group on connections

$$A\mapsto GAG^{-1}+\frac{dG}{dz}G^{-1}$$

\end{theorem}

The term "naturally isomorphic" means that this isomorphism preserver the action of the loop group.

Now take the one-dimensioanl $n$-th order scalar differential operator $L$.  As usual in theory of linear ODE to $L$ there corresponds a first-order multidimensional differential operator

\begin{equation}
\label{od}
d+\begin{pmatrix} 0 & u_2 &  u_3& ...& u_{n-1} & u_n\\
                              -1& 0 & 0 &...&0 &0\\
                              0&-1&0&...& 0&0\\
          ...\\
0&0&0&...&-1&0
\end{pmatrix}
\end{equation}

We can write it shortly as $d+J_{-}+A$, where

$$
J_{-}=\begin{pmatrix} 0 & 0 & 0& ...& 0 & 0\\
                              -1& 0 & 0 &...&0 &0\\
                              0&-1&0&...& 0&0\\
          ...\\
0&0&0&...&-1&0
\end{pmatrix},\,\,\,\,  A(x)=\begin{pmatrix} 0 & u_2 &  ...& u_{n-1} & u_n\\
                              0& 0 &...&0 &0\\
                              0&0&...& 0&0\\
          ...\\
0&0&...&0&0
\end{pmatrix}
$$

Now take a differential operator of the first order of type \eqref{do} and consider its orbit under the action of the loop group $GL(z)$. The orbit consists  of all operators with the same monodromy in $0$.

 Consider an action of the smaller group $N(z)$ defined in the previous section $GL(z)$-orbit decomposes into $N(z)$-orbits.  Every such orbit has a unique representative of type  \eqref{od} (see \cite{edfr1}).
Thus the set of operators of type  \eqref{od} is indeed a Poisson reduction of  the  dual to an affine Lie algebra.

\section{Poisson reduction: an approach based on constraints}

In this section we present another viewpoint to Poisson reduction  (see \cite{bai}). We describe this approach explicitly in our main example.


\subsection{Constraints of the first class}
Consider a Poisson manifold $M$ and some constraints of the first class  $\{\varphi_i\}$, $i=1,...,n$.
 These are some functions on $M$, such that

$$
\{\varphi_i,\varphi_j\}=0 \text{ on the set defined by equations }\varphi_i=0,\,\,i=1,...,n
$$

Equivalently, on the whole manifold $M$ one has

$$
\{\varphi_i,\varphi_j\}=c_{i,j}^k\varphi_k
$$

If one has a Poisson manifold with constraints $\varphi_i$ of the first type then one has an action of some Lie group on $M$. This is an exponential
 Lie group whose Lie algebra has generators  $x_i$ subject to relation $$[x_i,x_j]=c_{i,j}^k x_k,$$

the infinitesimal action of this group is defined by the following vector fields $X_i$ on $M$

\begin{equation}
\label{act}X_i=\{\varphi_i,.\}.
\end{equation}

This group is called {\it the gauge group}.

This action preserves obviously  the submanifold defined by constraints
$$
\{\varphi_i=0\}
$$

\subsection{Constraints associated with an $\mathfrak{sl}_2$-embedding}

The Lie algebra  $\mathfrak{sl}_2$ is generated as a linear space by three elements $e,h,f$ subject to relations

$$
[e,f]=h,\,\,\,[h,e]=2e,\,\,\,[h,j]=-2f.
$$

Consider an embedding

\begin{equation}
\label{ei}
i: \mathfrak{sl}_2\rightarrow  \mathfrak{g}
\end{equation}

Introduce a notation for an eigenspace of the operator $ad_h=[h,.]$

\begin{equation}
\mathfrak{g}_k=\{g\in \mathfrak{g},\,\,\, [h,g]=kg \},
\end{equation}

The eigenvalue $k$ is necessarily  an integer number.

 Then  for some $N\in \mathbb{Z}$ we have

\begin{equation}
\mathfrak{g}=\oplus_{k=-N}^N\mathfrak{g}_k
\end{equation}

Also an  embedding \eqref{ei} turns $\mathfrak{g}$ into a  $\mathfrak{sl}_2$-representation.
 This representation is reducible and can be decomposed into direct sum of irreducible representations.

Take a base in  $\mathfrak{g}$ denoted as  $I_{k,\mu,m}$, where $k$ is a highest weight of a $\mathfrak{sl}_2$-irreducible representation that contains the considered element, $m$ is a  $\mathfrak{sl}_2$-weight and the index $\mu$ indexes different $\mathfrak{sl}_2$-irreducible representation with the same highest weight.

As the underlying simple Lie algebra the affine Lie algebra $\widehat{\mathfrak{g}}$ has an invariant non-degenerate bilinear form (but it is not positively-defined, see \cite{vk}), denote it as $B$. We can suggest that the restriction of this form to $I_{\alpha}^0$ equals to the Killing form of $\mathfrak{g}$.
 Using the form $B$ we can identify $\widehat{\mathfrak{g}}$ and  $\widehat{\mathfrak{g}}^*$ according to the ruler

$$
x\leftrightarrow B(x,.).
$$
 The functions on $\widehat{\mathfrak{g}}^*$ then become functions on $\widehat{\mathfrak{g}}$, in particular the constraints are functions on $\widehat{\mathfrak{g}}$.

To define the constraints consider the loop algebra  $\mathfrak{g}((z))$. Every loop $l(z)\in\mathfrak{g}((t))$ can be written

\begin{equation}
l(z)=\sum_{k,m}U^{k,\mu,m}(z)I_{k,\mu,m},
\end{equation}

the coefficient  $U^{k,\mu,m}(z)=\sum_nU^{k,\mu,m}_nz^{-n-1} $ is a formal power series whose coefficients are complex numbers that depend linearly on $l$.  Let $I^{k,\mu,m}$ be a base in $\mathfrak{g}$ dual to  $I_{k,\mu,m}$. Then $U^{k,\mu,m}_n=I^{k,\mu,m}_n(l)$. We write also $U^{k,\mu,m}(z)=(I^{k,\mu,m}(z))(l)$.


Suppose that the image of $\mathfrak{sl}_2$   under the embedding has the highest weight $1$ and the index $\mu=1$, thus $I_{1,1,1}=t_+$, $I_{1,1,0}=t_0$, $I_{1,1,-1}=t_-$. Now take constraints


\begin{align}
\begin{split}
\label{constr}
\varphi^{k,\mu,m}(z)=I^{k,\mu,m}(z)-\delta_1^j\delta_1^m\delta^{\mu}_1,\,\,\, m\geq 0,
\end{split}
\end{align}

Here we mean that $\varphi^{k,\mu,m}(z)=\varphi^{k,\mu,m}_n z^{-n-1}$, and we take as constraints all functions  $\varphi^{k,\mu,m}_n $.
We identify  $I^{k,\mu,m}_n$ with a function $B(I^{k,\mu,m}_n,.)$.

\begin{prop}
These constraints are the first class
\end{prop}


\subsection{The group defined by constraints}

Denote as $\varphi^{\alpha}(y)$ a constraint defined in the formula \eqref{constr}. According to the general scheme we identify the modes
 modes $\varphi^{\alpha}_n$ with generators of some Lie algebra and we  define the action of these generators on a function on $\widehat{\mathfrak{g}}^*_1$ by vector fields define in the general formula \eqref{act}.

Using the bilinear form we identify  $\widehat{\mathfrak{g}}^*_1$ and $\widehat{\mathfrak{g}}_1$. We need to defined the action of $\varphi^{\alpha}_n$ onto an element of  $\widehat{\mathfrak{g}}_1$. This action must preserve the central element $1$, thus we need to define the action onto the loop $l(z)=\sum_{\beta}(I^{\beta}(z))(l)I_{\beta}$. Obviously we can take just $I^{\beta}(x)I_{\beta}$. Denote the result of this  infinitesimal action as $\delta_{\varphi^{\alpha}_n}I^{\beta}(x)I_{\beta}$

We can consider all $\varphi^{\alpha}_n$ simultaneously by taking $$\delta_{\varphi^{\alpha}(y)}I^{\beta}(x)I_{\beta}:=\sum_n y^{-n-1}\delta_{\varphi^{\alpha}_n}I^{\beta}(x)I_{\beta}$$

We have

\begin{align*}
& \delta_{\varphi^{\alpha}(y)}(I^{\beta}(x)I_{\beta})=\{\varphi^{\alpha}(y),  I^{\beta}(x)I_{\beta}\}=\{I^{\alpha}(y),  I^{\beta}(x)\}I_{\beta}=
\\&=B^{\alpha,\beta}\delta'(x-y)I_{\beta}+f^{\alpha,\beta}_{\gamma}\delta(x-y)I^{\gamma}(x)I_{\beta}.
\end{align*}

Here $B^{\alpha,\beta}=B(I_{\alpha},I_{\beta})$.
But since $f_{\alpha,\beta,\gamma}$ is antisymmetric one has

\begin{align*}
&f^{\alpha,\beta}_{\gamma}I_{\beta}=[I_{\gamma},I_{\delta}]B^{\alpha,\delta}.
\end{align*}

Thus one obtains

\begin{align*}
& \delta_{\varphi^{\alpha}(y)}(I^{\beta}(x)I_{\beta})=B^{\alpha,\beta}\delta'(x-y)I_{\beta}+[I_{\gamma},I_{\delta}]B^{\delta,\alpha}\delta(x-y)I^{\gamma}(x)
\end{align*}

Note that

 $B^{a,\alpha}I_{\alpha}\in \mathfrak{g}_{-k}$ if and only if $I_{\alpha}\in \mathfrak{g}_{k}$

The  conclusion is the following:

\begin{prop}
That is the Lie algebra of gauge transformations is

$$\oplus_{k\geq 1} \mathfrak{g}_{-k}.$$
\end{prop}

\begin{prop}
In the case $\mathfrak{g}=\mathfrak{gl}_N$, $t_{+}=E_{1,2}+E_{2,3}+...+E_{N-1,N}$ (then $t_{-}$ and $t_0$ are uniquely defined, see Section \ref{jmt}) the corresponding  group is  $N(z)$  is the group  of unipotent upper-triangular matrices whose coefficients depend on $z$.
\end{prop}

\begin{theorem}
This approach is equivalent to the one described in the previous section.
\end{theorem}

\section{Classical  finite $W$ algebras}

 The construction outlined in the previous section is the following:  a classical $W$ algebra is a Poisson algebra that is obtained from the Kirillov-Poisson structure on the dual to the affine Lie algebra using the Hamiltonian reduction.

There naturally rises  a question: shall we obtain something interesting if we take in this construction a simpler object: a dual space to a simple Lie algebra?

The obtained object is indeed interesting and it is called {\it a classical finite $W$ algebra}. This program was initiated in \cite{t}, in that paper firstly the classical finite $W$-algebras were defined (see also a review \cite{t1}).

\subsection{A description of the classical $W$-algebra in the second approach}

\label{clw}
Let us give a description of a classical $W$-algebra. We follow the paper \cite{bt}.

Take a simple Lie algebra $\mathfrak{g}$ and  fix an embedding

$$i:\mathfrak{sl}_2\rightarrow \mathfrak{g}$$

Denote the images of the elements $e,f,h$ as $t_{+},t_{-},t_0$.

Take a base in  $\mathfrak{g}$ denoted as  $I_{k,\mu,m}$, where $k$ is a highest weight of a $\mathfrak{sl}_2$-irreducible representation that contains the considered element, $m$ is a  $\mathfrak{sl}_2$-weight and the index $\mu$ indexes different $\mathfrak{sl}_2$-irreducible representation with the same highest weight.

The dual space $\mathfrak{g}^*$ has a structure of a Poisson manifold.
Using the Killing form $B(.,.)$ we can identify $\mathfrak{g}$ and $\mathfrak{g}^*$ by formula
$$
x\leftrightarrow B(x,.).
$$

Take in $\mathfrak{g}$ a base $I^{k,\mu,m}$ dual to $I_{k,\mu,m}$.
One has

\begin{align}
\{I^{\alpha},I^{\beta}\}=f_{\gamma}^{\alpha,\beta}I^{\gamma}
\end{align}

The functions on $\mathfrak{g}^*$ then become functions on $\mathfrak{g}$, in particular the constraints are functions on $\mathfrak{g}$.
Consider the following constraints

\begin{equation}
\label{constrfin}
\varphi_{k,m}^{\mu}=I^{k,\mu,m}-\delta_{1}^k\delta_{1}^m\delta_{\mu}^1,\,\,\,\, m>0
\end{equation}
 Here  $I^{k,\mu,m}$ is considered as a function $B(I^{k,\mu,m},.)$ on  $\mathfrak{g}$.

Consider the subset $\mathfrak{g}_c$ in $\mathfrak{g}$ defined by equations

$$
\varphi_{k,m}^{\mu}=0
$$

Explicitly the elements of this set are given by the formula

$$
t_++\sum_{k,\mu}\sum_{m<0}\alpha_{j,m}^{\mu}I_{j,m}^{\mu},
$$

where $\alpha_{j,m}^{\mu}$ are arbitrary constants.

\begin{prop}
These constant are of the first class.
\end{prop}

Now let us  write explicitly the group action generated by these constraints. This is an action on the function on $\mathfrak{g}^*$. A function on $\mathfrak{g}^*$ is an element of $\mathfrak{g}^{**}=\mathfrak{g}$.

Let $\varphi^{\alpha}$ be a constraint defined in the formula \eqref{constrfin}.   According to the general scheme we identify the modes
 modes $\varphi^{\alpha}$ with generators of some Lie algebra and we  we  define the action of these generators on a function on $\mathfrak{g}$ by vector fields define in the general formula \eqref{act}.  We denote the result of the action  onto a function $x$ as $\delta_{\varphi^{\alpha}} x$. It is sufficient to defined an action of  $\varphi^{\alpha}$
onto an element  $x$ proportional to $I_{\beta}$. Such an element can be written as $x=x^{\beta}I_{\beta}\in \mathfrak{g}$. And $x^{\beta}$ can be written as $I^{\beta}(x)$, i.e. the value of the function $I^{\beta}$ on the elements $x$.

 One has

\begin{align*}
&\delta_{\varphi^{\alpha}} x=\{\varphi^{\alpha},x^{\beta}I_{\beta}\}=\{I^{\alpha},I^{\beta}(x)\}I_{\beta}=f^{\alpha,\beta}_{\gamma}I^{\gamma}I_{\beta}
\end{align*}

But since $f_{\alpha,\beta,\gamma}$ is antisymmetric one has

\begin{align*}
&f^{\alpha,\beta}_{\gamma}I_{\beta}=[I_{\gamma},I_{\delta}]B^{\alpha,\delta},
\end{align*}

where $B^{\alpha,\delta}=B(I_{\alpha},I_{\beta})$.
Thus one obtains

\begin{align*}
&\delta_{\varphi^{\alpha}} x=[I_{\gamma},I_{\delta}]B^{\alpha,\delta}J^{\gamma}(x)
\end{align*}

Since $B^{\alpha,\delta}I_{\delta}\in \mathfrak{g}^{-k}$ if and only if $I_c\in\mathfrak{g}^{k}$ one find that the Lie algebra of gauge transformations is

$$\oplus_{k\geq 1}\mathfrak{g}^{k}.$$

\begin{prop}
In the case $\mathfrak{g}=\mathfrak{gl}_N$, $t_{+}=E_{1,2}+E_{2,3}+...+E_{N-1,N}$ (then $t_{-}$ and $t_0$ are uniquely defined, see Section \ref{jmt}) the group  $N$  of unipotent upper-triangular matrices.
\end{prop}

The Poisson reduction is the factor space  $\mathfrak{g}_c/N$.

\begin{prop}
Every $N$-orbit has a unique representative of type

$$
t_{+}=\sum_{j,\mu}x_{j}^{\mu}I_{j,-j}^{\mu},\,\,\,\,\,x_{j}^{\mu}\in\mathbb{C}
$$
\end{prop}

Thus the factor space is isomorphic to the space

\begin{equation}
\label{glw}
\mathfrak{g}_{lw}=<  \sum_{j,\mu}x_{j}^{\mu}I_{j,-j}^{\mu} >
\end{equation}

which is just the space formed by $\mathfrak{sl}_2$-lowest vectors.

Let us give a description of a Poisson bracket on this space.

To do it let us  define an operator $L$.  We have a map

$$
ad_{t_+}: Im(ad_{t_{-}})\rightarrow Im (ad_{t_{+}})
$$

Let us denote the inverse operator continued to the rest part of $\mathfrak{g}$ by $0$ as $L$.

Then for functions $I^{\alpha}, I^{\beta}\in C^{\infty}(\mathfrak{g}_{lw})$ one has

\begin{equation}
\label{scopp}
\{I^{\alpha},I^{\beta}\}(I_{\gamma})=B(I_{\gamma},[I^{\alpha},\frac{1}{1+L\circ ad_{I_{\gamma}}}I^{\beta}])
\end{equation}

Here $I^{\alpha}$ is an element of  a base of $\mathfrak{g}$ dual to the base $I_{\alpha}$.

Thus we come to the following statement

\begin{theorem}
The classical $W$ algebra associated to the embedding $i$ is the space
\begin{equation}
\mathfrak{g}_{lw}<   \sum_{j,\mu}x_{j}^{\mu}I_{j,-j}^{\mu}>
\end{equation}
with the Poisson structure

\begin{equation}
\{Q_1,Q_2\}(w)=(w,[grad_wQ_1,\frac{1}{1+L\circ ad_w}grad_wQ_2]),
\end{equation}

where $Q_1(w)$, $Q_2(w)$ are arbitrary function on $\mathfrak{g}_{lw}$ and $w$  is a coordinate on  $\mathfrak{g}_{lw}$.
\end{theorem}

Consider a trivial example. If $t_{+}=0$ then all vectors are $\mathfrak{sl}_2$-lowest. Also $L=0$ and we obtain than $\mathfrak{g}_{lw}=\mathfrak{g}$ and $L=0$. Hence

$$
\{I^{\alpha},I^{\beta}\}(I_{\gamma})=B(I_{\gamma},[I^{\alpha},I^{\beta}])=f^{\alpha,\beta}_{\gamma}
$$

Thus the corresponding $W$-algebra is just $\mathfrak{g}$ with a canonical Poisson structure.

\subsection{Examples}

Take $\mathfrak{g}=\mathfrak{sl}_2=<f,h,e>$ and $i=id$.

Then $\mathfrak{g}_{lw}=<f>$.  Since the Poisson bracket is anti-symmetric we obtain that
the classical $W$-algebra corresponding to $\mathfrak{sl}_2$ and a trivial embedding of $\mathfrak{sl}_2$ into itself is a commutative one-dimensional Poisson algebra.

The simplest non-principle embedding is an embedding $i:\mathfrak{sl}_2\rightarrow \mathfrak{sl}_3$ defined by formulas

\begin{align*}
& t_{+}=E_{1,3},\,\,\,\, t_0=diag(\frac{1}{2},0,-\frac{1}{2}),\,\,\,\, t_{-}=E_{3,1}
\end{align*}

In this case the constraints are

\begin{align*}
&\varphi_{1,1}^1=I_{1,1}^1-1=0\\
&\varphi_{1/2,1/2}^1=I_{1/2,1/2}^1=0,\\
&\varphi_{1/2,1/2}^2=I_{1/2,1/2}^2=0.
\end{align*}

Choose the following generators of $(\mathfrak{sl}_3)_{lw}$

\begin{align*}
&\mathfrak{g}_{lw}=<\frac{1}{6}(E_{1,1}+E_{3,3}-2E_{2,2}),E_{1,2},E_{2,3},E_{1,3}>
\end{align*}

According to our notations we denote these elements as

$$
I^1_{0,0}, I^1_{1/2,-1/2}, I^2_{1/2,-1/2}, I^1_{1,-1}
$$

\begin{align*}
&C=-\frac{4}{3}(I_{1,-1}^1+3(I_{0,0}^1)^2,\\
&E=I_{1/2,1/2}^1,\\
&F=\frac{4}{3}I^2_{1/2,-1/2},\\
&H=4I_{0,0}^1
\end{align*}

Then the Poisson brackets are

\begin{align*}
&\{H,E\}=2E,\\
&\{H,F\}=-2F\\
&\{E,F\}=H^2+C\\
&C\text{ commutes with everything }
\end{align*}

An embedding $i:\mathfrak{sl}_2\rightarrow \mathfrak{g}$ is called principal if it's coadjoit orbit in $\mathfrak{g}$ is of maximal dimension.
In the case  of $\mathfrak{g}=\mathfrak{gl}_{n}$ this means that $e$ is conjugate to

$$
\begin{pmatrix}
0& 1 &0 &...&0\\
0&0& 1& ...&0\\
...\\
0&0&0&...&1\\
0&0&0&...&0
\end{pmatrix}
$$

\begin{prop}
A classical $W$-algebra corresponding to a principal embedding is a commutative Poisson
algebra.
\end{prop}

\section{A definition of a finite $W$-algebra}
Above we have define a classical $W$ algebras. What is the construction of it's quantum analogue?

One can prove that the quantization of the Poisson algebra of functions on $\widehat{\mathfrak{g}}^*_1$ is the algebra $\widehat{\mathfrak{g}}_1$. But we have a reduction of that Poisson structure defined by constraints of the first type. The procedure of quantization of a Poisson algebra which is obtained by imposing constraints of the first type is known in physics - this is  the BRST procedure.
If one applies it to infinite classical $W$ algebra one obtains an infinite $W$ algebra \cite{ff}.

The first definition of the finite $W$ algebra were given in this way \cite{bt}. So this is the paper, where finite $W$-algebras were discovered.
But later there were obtained simpler definitions. To explain then we need some facts bout $\mathfrak{sl}_2$-triples and nilpotent orbits.

\subsection{Nilpotent orbits}

Take our  Lie algebra $\mathfrak{g}$. An element $e\in\mathfrak{g}$ is nilpotent if $ad_{e}$
is a nilpotent endomorphism of $\mathfrak{g}$.For example in $\mathfrak{gl}_n$ the elements $E_{i,j}$, $i\neq j$ are nilpotent.

For a arbitrary $x\in\mathfrak{g}$ denote as $\mathfrak{g}^x$ the kernel of the mapping $ad_x$, in other words this is a centralizer of $x$. This is a subalgebra in $\mathfrak{g}$.

\subsection{$\mathbb{Z}$-grading}

A $\mathbb{Z}$-grading of a Lie algebra $\mathfrak{g}$ is a decomposition into a direct sum

\begin{align*}
&\mathfrak{g}=\oplus_{j\in\mathbb{Z}}\mathfrak{g}_j\text{ such that }\\
&[\mathfrak{g}_i,\mathfrak{g}_j]\subset\mathfrak{g}_{i+j}.
\end{align*}

Denote as $G$ an exponential Lie group corresponding to $\mathfrak{g}$.  One has an action of $G$ on $\mathfrak{g}$, which is described  by the formula

$$x\mapsto  G^{-1}xG,$$

denote an orbit of $e$ as $\mathcal{O}_e$, this a subspace in $\mathfrak{g}$.

There exists a unique dense nilpotent  orbit called {\it the regular nilpotent orbit}, the corresponding nilpotent element is called regular nilpotent.
Equivalently one ca say that $e$ is regular nilpotent if $\mathfrak{g}^e$ is  of minimal dimension.

The set of all nilpotent orbits is  a partially ordered set

\begin{equation}
\mathcal{O}\leq \mathcal{O}' \leftrightarrow \mathcal{O}\subset \text{closure of }\mathcal{O}'
\end{equation}

For example in $\mathfrak{gl}_N$ the nilpotent orbit are parameterized by partitions $\lambda=(\lambda_1\geq \lambda_2\geq...)$, $N=\lambda_1+\lambda_2$ of $N$. A partition corresponds to an orbit that contains a nilpotent  matrix in Jordan form $(J_{\lambda_1},...)$. One has $\mathcal{O}_{reg}\leftrightarrow (N)$, $\mathcal{O}_{sub}\leftrightarrow (N,1)$ is a unique dense orbit in $\mathfrak{gl}_N\setminus \mathcal{O}_{reg}$.

\subsection{Jacobson-Morozov theorem}
\label{jmt}
The theorem of Jacobson-Morozov says the  following

\begin{theorem} Every non-zero nilpotent element $e$ can be included in a $\mathfrak{sl}_2$-triple $\{e,f,h\}$. That is elements these elements satisfy the  $\mathfrak{sl}_2$- commutation relations

$$
[e,f]=h,\,\,\,[h,e]=2e,\,\,\,[h,j]=-2f.
$$

\end{theorem}

For example if $\mathfrak{g}=\mathfrak{gl}_{N}$, $e=J_N$ - a regular nilpotent element, then

\begin{align*}
&h=diag(N-1,N-3,...,3-N,1-N),\\
&f=\begin{pmatrix}  0 & 0 &....& 0 & 0\\  a_1 & 0& ...& 0 &0 \\ ...\\  0&0&...&a_{N-1} &0 \end{pmatrix},\,\,\, a_i=i(n-i)
\end{align*}

\subsection{Gragings}

If we are given an  $\mathfrak{sl}_2$- triple $\{e,f,h\}$, then we have a decomposition

\begin{equation}
\mathfrak{g}=\oplus\mathfrak{g}_i,\,\,\,\mathfrak{g}_i=\{g\in\mathfrak{g}:\,\, [h,g]=jg\},
\end{equation}

this grading is called {\it a Dynkin grading}, it satisfies the following properties

\begin{align}
& e\in\mathfrak{g}_2,\\
&ad_e: \mathfrak{g}_{j}\rightarrow\mathfrak{g}_{j+2} \text{ is injective for }j\leq -1,\\
&ad_e : \mathfrak{g}_{j}\rightarrow\mathfrak{g}_{j+2} \text{ is surjective for }j\geq -1,\\
&\mathfrak{g}_e\subset \oplus_{j\geq 0}\mathfrak{g}_j,\\
&B(\mathfrak{g}_j,\mathfrak{g}_i)=0\text{ unless }i+j=0,\\
& dim\mathfrak{g}_e=dim\mathfrak{g}_0+dim\mathfrak{g}_1.
\end{align}

An arbitrary grading is called {\it good} if it satisfies properties

\begin{align}
& e\in\mathfrak{g}_2,\\
&ad_e: \mathfrak{g}_{j}\rightarrow\mathfrak{g}_{j+2} \text{ is injective for }j\leq -1,\\
&ad_e : \mathfrak{g}_{j}\rightarrow\mathfrak{g}_{j+2} \text{ is surjective for }j\geq -1,\\
&\mathfrak{g}_e\subset \oplus_{j\geq 0}\mathfrak{g}_j,\\
&B(\mathfrak{g}_j,\mathfrak{g}_i)=0\text{ unless }i+j=0,\\
& dim\mathfrak{g}_e=dim\mathfrak{g}_0+dim\mathfrak{g}_1.
\end{align}

i.e. (1)-(6) above.

For example take $\mathfrak{g}=\mathfrak{gl}_3$, $e=E_{1,3}$.

\begin{enumerate}
\item Let $h=diag(1,0,-1)$, $f=E_{3,1}$, then the degrees of $E_{i,j}$ in the Dynkin grading can be presented as follows

$$
\begin{pmatrix}
0&1&2\\
-1&0&1\\
-2&-1&0
\end{pmatrix}
$$

\item Take an element $h=diag(1,1,-1)$ and consider it's eigenspace decomposition, then we obtain a good, but non-Dynkin grading, defined by the matrix

$$
\begin{pmatrix}
0&0&2\\
0&0&2\\
-2&-2&0
\end{pmatrix}
$$

\end{enumerate}

\subsection{A bijection between nilpotent orbits and  $\mathfrak{sl}_2$-triples}

\begin{theorem}The mapping

\begin{align*}
& \{\mathfrak{sl}_2-\text{ triples}\}/G\rightarrow \{\text{non-zero nilpotent orbts}\}\\
&\{e,f,h\}\mapsto \mathcal{O}_e
\end{align*}
is a bijection
\end{theorem}

\subsection{Definition  of finite $W$ algebras via Whittaker modules}

\subsubsection{A definition}
Now we can give the first definition of a finite $W$ algebra.
This definition was given in \cite{p1}

So take a reductive Lie algebra  $\mathfrak{g}$, an nilpotent element $e$. Take an algebra

$$\mathfrak{sl}_2=<e,f,h>$$
and an embedding
\begin{equation}
i:\mathfrak{sl}_2 \hookrightarrow \mathfrak{g}
\end{equation}

Consider a linear form

\begin{equation}
\chi(x):=B(x,e).
\end{equation}

Introduce a notation

\begin{align*}
&<.,.>\mathfrak{g}_{-1}\times  \mathfrak{g}_{-1}\rightarrow \mathbb{C}\\
&<x,y>:=B([x,y],e)=\chi([x,y])
\end{align*}

\begin{prop}
This form is skew-symmetric and non-degenerate
\end{prop}

Denote as $\mathfrak{l}$ the Lagrangian subspace in $\mathfrak{g}_{-1}$ that is the maximal
subspace such that $<.,.>\mid_{\mathfrak{l}}=0$. One has $dim\mathfrak{l}=\frac{1}{2}dim\mathfrak{g}_{-1}$.

Put

\begin{equation}
\mathfrak{m}=\mathfrak{l}\oplus\bigoplus_{j\leq-2}\mathfrak{g}_{j}
\end{equation}

It is a Lie subalgebra in $\mathfrak{g}$.

\begin{prop}
The mapping $$\chi\mid_{\mathfrak{m}}:\mathfrak{m}\rightarrow \mathbb{C}$$
is a one-dimensional representation of $\mathfrak{m}$
\end{prop}

This Proposition is an immediate consequence of  the fact that $\mathfrak{l}$ is lagrangian.

Denote as $I_{\chi}$ the following left ideal in $U(\mathfrak{g})$:

\begin{align*}
&I_{\chi}:=\text{ left ideal generated by }a-\chi(a),\,\,\, a\in\mathfrak{m}=\\
&=\{\ x(a-\chi(a)),\,\,\,x\in U(\mathfrak{g})\,\,\, a\in\mathfrak{m}\}
\end{align*}

Since $\mathfrak{m}$ act on $\mathfrak{g}$ in an  adjoint way, $U(\mathfrak{g})$ is a $U(\mathfrak{m})$-module. Consider an induced module
\begin{equation}
Q_{\chi}:=U(\mathfrak{g})\otimes_{U(\mathfrak{m})}\mathbb{C}=U(\mathfrak{g})/I_{\chi}
\end{equation}

A finite $W$-algebra is defined as follows

\begin{equation}
W_{\chi}:=End_{U(\mathfrak{g})}(Q_{\chi})^{op}
\end{equation}

Consider a trivial example.  Let $e=0$, then $\chi=0$, $\mathfrak{g}=\mathfrak{g}_0$, hence $\mathfrak{m}=0$ and $Q_{\chi}=U(\mathfrak{g})$ and $W_{\chi}=U(\mathfrak{g})$.  Compare this example with a trivial example at the end of section \ref{clw}.

\subsubsection{A  reformulation of the definition}
\label{sedef}
Note that the $U(\mathfrak{g})$-module $Q_{\chi}=U(\mathfrak{g})/I_{\chi}$ is cyclic, hence it's endomorphism is defined by an image of the coset $\bar{1}:=1+I_{\chi}$, and this image must be annihilated by $I_{\chi}$. So

$$
W_{\chi}=\{\bar{y}\in U(\mathfrak{g})/I_{\chi}:\,\,\, (a-\chi(a))y\in I_{\chi},\,\,\,\forall a\in\mathfrak{m}\}
$$

We can reformulate this as follows

\begin{equation}
W_{\chi}=(Q_{\chi})^{ad\mathfrak{m}}=\{\bar{y}\in U(\mathfrak{g})/I_{\chi}:\,\,\, [a,y]\in I_{\chi}\,\,\, \forall a\in \mathfrak{m}\}
\end{equation}

The algebra structure is given by

$$
\bar{y}_1\bar{y}_2:=\overline{y_1y_2}
$$

\subsection{A second definition of a finite $W$-algebra}

\subsubsection{A definition}

This definition was given in  \cite{ko} even before the finite $W$-algebras were discovered.
Consider the case when the grading is even, that is $\mathfrak{g}=0$ unless $j$ is even. {\bf Everywhere below in the paper we suggest that the grading is even.}

Then $\mathfrak{g}_{-1}=0$, hence $\mathfrak{m}=\bigoplus_{j\leq -2}\mathfrak{g}_{j}$. Put

$$
\mathfrak{v}=\bigoplus_{j\geq 0}\mathfrak{g}_j
$$

This is a subalgebra in $\mathfrak{g}$.
One has

$$
U(\mathfrak{g})=U(\mathfrak{v})\oplus I_{\chi}
$$

Denote as $pr_{\chi}$ the projection to the first summand. It defines an isomorphism

\begin{equation}
\label{wm2}
\overline{pr_{\chi}}: U(\mathfrak{g})/I_{\chi}\rightarrow U(\mathfrak{v}).
\end{equation}

The first definition of a finite $W$ algebra was

\begin{equation}
\label{wm1}
W_{\chi}=(Q_{\chi})^{ad\mathfrak{m}}=\{\bar{y}\in U(\mathfrak{g})/I_{\chi}:\,\,\, [a,y]\in I_{\chi}\,\,\, \forall a\in \mathfrak{m}\}
\end{equation}

Comparing \eqref{wm1} and \eqref{wm2} we come to the definition

\begin{equation}
\label{wv}
W_{\chi}=U(\mathfrak{v})^{ad\mathfrak{m}}=\{y\in U(\mathfrak{v}):\,\,\, [a,y]\in I_{\chi} \forall a\in\mathfrak{m}\}
\end{equation}

\subsubsection{An example}

Take $e=E_{1,2}\in \mathfrak{gl}_2$, then $\mathfrak{m}=\mathbb{C}f$, where $f=E_{2,1}$.

 Also $\mathfrak{v}=\mathbb{C}e+\mathbb{C}E_{1,1}+\mathbb{C}E_{2,2}$.

One can show that $W_{\chi}$ is  a polynomial algebra generated by $E_{1,1}+E_{2,2}$ and $e+\frac{1}{4}h^2-\frac{1}{2}h$.

The simplest non-principle embedding is an embedding $i:\mathfrak{sl}_2\rightarrow \mathfrak{sl}_3$ defined by formulas

\begin{align*}
& t_{+}=E_{1,3},\,\,\,\, t_0=diag(\frac{1}{2},0,-\frac{1}{2}),\,\,\,\, t_{-}=E_{3,1}
\end{align*}

The degrees of elements of $\mathfrak{sl}_3$ are given in the matrix

$$
\begin{pmatrix}
0&1&-2\\
-1&0&1\\
-2&-1&0
\end{pmatrix}
$$

Thus

\begin{align*}
&\mathfrak{g}^e=<\frac{1}{6}(E_{1,1}+E_{3,3}-2E_{2,2}),\frac{1}{2}(E_{3,3}-E_{1,1}),E_{1,2},E_{2,3},E_{1,3}>
\end{align*}

Denote these elements as

$$
I^1_{0,0},I^1_{1,0}, I^1_{1/2,-1/2}, I^2_{1/2,-1/2}, I^1_{1,-1}
$$

Choose the following generators of $(\mathfrak{sl}_3)_{lw}$

\begin{align*}
&C=-\frac{4}{3}(I_{1,-1}^1+3(I_{0,0}^1)^2),\\
&E=I_{1/2,-1/2}^1,\\
&F=\frac{4}{3}I^2_{1/2,-1/2},\\
&H=4I_{0,0}^1
\end{align*}

Then the Poisson brackets are

\begin{align*}
&[H,E]=2E,\\
&[H,F]=-2F\\
&[E,F]=H^2+C\\
&C\text{ commutes with everything }
\end{align*}

In Section \ref{cent} we present a description of a $W$-algebra associated with a principle nilpotent.

\subsection{Definition though the BRST procedure}

This definition makes a bridge to the previously  discussed classical finite $W$-algebras.

Consider for simplicity only the case of an even grading.  Then $\mathfrak{m}=\oplus_{j\leq -2}\mathfrak{g}_{j}$.
Take a copy of $\mathfrak{m}$ and denote it as $\hat{\mathfrak{m}}$. Also take a dual space $\mathfrak{m}^*$.

On the direct sum $\mathfrak{m}^*\oplus \hat{\mathfrak{m}}$ with a  symmetric bilinear form defined by pairing. The Clifford algebra of this space is $\Lambda(\mathfrak{m}^*)\otimes \Lambda( \hat{\mathfrak{m}})$. And the tensor algebra of this Clifford algebra is

\begin{equation}
\mathcal{B}^*=\Lambda(\mathfrak{m}^*)\otimes U(\mathfrak{g}) \otimes  \Lambda( \hat{\mathfrak{m}})
\end{equation}


Take a basis $\{b_i\}$ in $\hat{\mathfrak{m}}$ and a dual basis $\{f_i\}$
 in $\mathfrak{m}^*$. Introduce an element

\begin{equation}
\varphi=f^i(b_i-\chi(b_i))-\frac{1}{2}f^if^j[b_i,b_j]^{*}\in  \mathcal{B}^*,
\end{equation}

where $a^{*}$ is a dual element.

Let $d=[\varphi,.]$.

Explicitely

\begin{align*}
&d(x)=f^i[b_i,x],\,\,\, x\in\mathfrak{g},\\
&d(f)=\frac{1}{2}f^iad^*b_i(f),\\
&d(a^*)=a-\chi(a)+f^i[b_i,a]^*
\end{align*}

Here $ad^*$ is a coadjoint action.

\begin{prop}
$d^2=0$
\end{prop}

\begin{theorem}
$$
H^0(\mathcal{B}^*)=W_{\chi}
$$
\end{theorem}

\section{The structure of a finite $W$-algebra}
In this section we present some result about
 the structure of a finite $W$-algebra. In particular we formulate the analogues of PBW theorem.

\subsection{The Kazhdan filtration}

We consider a finite $W$-algebra in a Whittaker-module realization.
A filtration of an algebra $A$ is a sequence of vector subspaces

$$
...\subset F_i A\subset F_{i+1} A\subset...
$$

such that

\begin{enumerate}
\item $\cap_i F^i A=0$, $\cup_i F_iA=A$
\item $F_jA\cdot F_iA\subset F_{i+j}A$
\end{enumerate}

The associated graded algebra is the algebra
$$
gr(A)=\oplus_i (F_iA/F_{i+1}A)
$$

The algebra $U(\mathfrak{g})$ has a filtration

$$...\subset F_iU(\mathfrak{g})\subset F_{i+1}U(\mathfrak{g})\subset...$$
 called {\it a Kazhdan filtration}, which is defined as follows.

An embedding $$i:\mathfrak{sl}_2\rightarrow \mathfrak{g}$$ induces a decomposition

$$\mathfrak{g}=\oplus \mathfrak{g}_i,\,\,\,\mathfrak{g}_i=\{g\in\mathfrak{g}:\,\,[h,g]=ig\}$$

 The elements $x\in\mathfrak{g}_i$ have degree $i+2$, and the products $x_1...x_n$, where $x_k\in\mathfrak{g}_{i_k}$ have degree
$(i_1+2)+...+(i_n+2)$.

Note that the associated graded $gr U(\mathfrak{g})$ is a symmetric algebra $S(\mathfrak{g})$ but with a non-standard grading:  $x\in\mathfrak{g}_i$ have degree $i+2$.

The Kazhdan filtration on $U(\mathfrak{g})$ induces a filtration  on $I_{\chi}$ and $I_{\chi}$ is a graded ideal in $ U(\mathfrak{g})$.  Hence we have an induced on $Q_{\chi}=U(\mathfrak{g})/I_{\chi}$ and on $W_{\chi}=End(Q_{\chi})^{op}$:

$$
F^iW_{\chi}=\{\ w\in W_{\chi}:\,\,\, w(F^kQ_{\chi})\subset F^{k+i} Q_{\chi} ) \}.
$$

\subsection{The associated graded to Kazhdan filtration}
 To give a description of the associated graded we use the definition

\begin{equation}
W_{\chi}=(Q_{\chi})^{ad\mathfrak{m}}=\{\bar{y}\in U(\mathfrak{g})/I_{\chi}:\,\,\, [a,y]\in I_{\chi}\,\,\, \forall a\in \mathfrak{m}\}
\end{equation}

and  make an identification

$$
grU(\mathfrak{g})=S(\mathfrak{g})=\text{ the  function on the affine variety of }\mathfrak{g}
$$

Also remind that

\begin{align*}
&I_{\chi}=\text{ left ideal in $U(\mathfrak{g})$ generated by }a-\chi(a),\,\,\, a\in\mathfrak{m}=\\
&=\{\ x(a-\chi(a)),\,\,\,x\in U(\mathfrak{g})\,\,\, a\in\mathfrak{m}\}
\end{align*}

The associated graded of $I_{\chi}$ is then described as an ideal in this algebra defined  by the following condition

$$
grI_{\chi}=\text{ ideal of functions   vanishing on the closed subvariety }e+\mathfrak{m}^{+},
$$

where $\mathfrak{m}^{+}$ is the orthogonal compliment of $\mathfrak{m}$ in $\mathfrak{g}$ with respect to the bilinear form $B(.,.)$.

Then $grQ=grU(\mathfrak{g})/ gr I_{\chi}$ is described as follows

$$
grQ=\text{ function on }e+\mathfrak{m}^{+}=\mathbb{C}[e+\mathfrak{m}^{+}]
$$

Note that from the second definition we obtain immediately

$$
grQ=S(\mathfrak{v}).
$$

But our purpose is to obtain a description of $grW_{\chi}=grQ^{\mathfrak{m}}$.

Take a subgroup $H$ in $G$ corresponding to a subalgebra $\mathfrak{m}$. It's adjoint action on $\mathfrak{g}$ leaves  $e+\mathfrak{m}^{+}$ invariant, so $H$ acts on  the space of functions $\mathbb{C}[e+\mathfrak{m}^{+}]$.

Put

$$
\mathfrak{g}^f=\{g\in \mathfrak{g}:\,\,\,\,[f,g]=0\}.
$$

We call the set $$e+\mathfrak{g}^f$$ {\it the  Slodowy slice.}

By  \cite{GG} the set $$e+\mathfrak{g}^f\subset e+\mathfrak{m}^+$$  gives a representative of each orbit. Hence

\begin{equation}
\label{mgf}
grQ^{\mathfrak{n}}=\mathbb{C}[e+\mathfrak{m}^+]^f\simeq\mathbb{C}[e+\mathfrak{g}^f].
\end{equation}

They say that the $W_{\chi}$ algebra is a quantization of the Slodowy slice.

\subsection{An analogue of a PBW theorem}

Take the algebra $U(\mathfrak{g})$ with a usual filtration

$$
F_i(U(\mathfrak{g}))=<x_1....x_i>,\,\,\,\, x_i\in\mathfrak{g}.
$$

Then $gr(U(\mathfrak{g}))=S(\mathfrak{g})$

The famous PBW theorem says that is we take a base $x_1,...,x_n$ in $\mathfrak{g}$ and fix some order $x_1<...<x_n$ then the products

$$
x_{i_1}...x_{i_k},\,\,\, x_{i_1}<...<x_{i_k}
$$

form a linear base of $F_{k}U(\mathfrak{g})$.

We can reformulate it as follows.

\begin{theorem}
There exist a  linear mapping

$$
\Theta: S(\mathfrak{g})=grU(\mathfrak{g})\rightarrow U(\mathfrak{g})
$$
such that  $\Theta (F_i(S(\mathfrak{g})))=  F_i(U(\mathfrak{g}))$.

This mapping commutes with the adjoint action of $\mathfrak{t}$.
\end{theorem}

One just puts

$$
\Theta(x_{i_1}...x_{i_k})=x_{(i_1)}...x_{(i_k)}, \text{ where } x_{(i_1)}<...<x_{(i_k)}
$$

Now let us turn to the algebra $W_{\chi}$. We have the following subalgebras in $\mathfrak{g}$:

$$
\mathfrak{g}^e=\{g\in \mathfrak{g}:\,\,\,\,[e,g]=0\},\,\,\mathfrak{t}^e=\{t\in \mathfrak{t}:\,\,\,\,[e,t]=0\},
$$

where  $ \mathfrak{t}$ is the Cartan subalgebra in $\mathfrak{g}$.

\begin{theorem}
There exists a $\mathfrak{t}^e$ equivarent map
$$
\Theta: \mathfrak{g}^e\rightarrow W_{\chi},
$$

such that for  $x\in \mathfrak{g}_j$ one has $\Theta(x)\in F_{j+2}(W)$ and such that if $x_1,...,x_t$ is a homogeneous base of $\mathfrak{g}^e$ (i.e. $x_i\in \mathfrak{g}^e\cap\mathfrak{g}_{n_i}$ for some $n_i\in\mathbb{Z}$) then

\begin{align*}
&\{\Theta(x_{i_1})...\Theta(x_{i_k}),\,\,\ k>0 ,\,\,\,1\leq i_1\leq...\\&...\leq i_k\leq t ,\,\,\, n_{i_1}+...+n_{i_k}+2k\leq j  \}
\end{align*}

form  a basis of $F_j W_{\chi}$.
\end{theorem}

Let us explain a construction of $\Theta$.  We have announced an isomorphism.
\begin{equation}
grW_{\chi}=\mathbb{C}[e+\mathfrak{m}^+]^f\simeq\mathbb{C}[e+\mathfrak{g}^f].
\end{equation}

Also the following direct sum decomposition takes place

\begin{lem}
$$
\mathfrak{v}=\mathfrak{g}^e\oplus\bigoplus_{j\geq 2}[f,\mathfrak{g}_j].
$$
\end{lem}

Take a projection to the first summand and consider a mapping
$$\zeta:S(\mathfrak{v})\rightarrow S(\mathfrak{g}^e) $$ induced by a projection to the first summand. Then since $B$ is invariant
for $x\in \bigoplus_{j\geq 2}\mathfrak{g}_j$ and $z\in e+\mathfrak{g}^f$ we have

\begin{equation}
\label{bxf}
B([f,x],z)=B(x,[z,f])=B(x,[e,f])=B(x,h)=0.
\end{equation}

We have used that $h\in\mathfrak{g}_0$, $x\in \bigoplus_{j\geq 2}\mathfrak{g}_j$ and these subspaces are orthogonal with respect to $B$.

The formula \eqref{bxf} means that an element from $ker\zeta$ annihilates $e+\mathfrak{g}^f$. Thus we obtain a mapping  from $Im\zeta=S(\mathfrak{g}^e) $ to $\mathbb{C}[e+\mathfrak{g}^f]$ which maps $z$ to $B(z,.)$.

\begin{lem}
This mapping is an isomorphism.
\end{lem}

Thus
we obtain that $grW_{\chi}=S(\mathfrak{g}^e)$,

 Introduce a notation for this isomorphism

\begin{equation}
\xi : grW_{\chi}\rightarrow  S(\mathfrak{g}^e)
\end{equation}

Now take $x_i\in \mathfrak{g}^e$, then we can take

$$
\Theta(x_i)=\xi^{-1}(x_i).
$$

\subsection{ A good filtration}
There exist another filtration on $W_{\chi}$. To define take a definition  given in the formula

\begin{equation}
W_{\chi}=U(\mathfrak{v})^{ad\mathfrak{m}}=\{y\in U(\mathfrak{v}):\,\,\, [a,y]\in I_{\chi} \forall a\in\mathfrak{m}\}
\end{equation}

Take a filtration  on $U(\mathfrak{v})$, and put

$$
F'_jW=W\cap F_jU(\mathfrak{v}).
$$

One has

\begin{prop}
$$
gr'(W)\subset U(\mathfrak{v})
$$

is a graded subalgebra
\end{prop}

Actually the following statement takes place

\begin{theorem}
$$
gr'(W)=U(\mathfrak{g}^e)
$$
\end{theorem}

\section{Representations of a finite $W$-algebra}

In this section we present results about representation of finite $W$ algebras. We present several independent approaches to  the subject.

\subsection{Primitive ideals}

Let us describe an approach to description of finite-dimensional irreducible representations of finite $W$-algebras that belongs to Losev \cite{lo}.
For a ring $A$ a primitive ideal is an ideal $I$ such that
 $$I=Ann(V)=\{a\in A:\,\,\, aV=0\}$$
  for some  irreduible $A$-module $V$.

Let us say some words about primite ideals in $U(\mathfrak{g})$. These ideals are classified (see \cite{Du}).  They have two main invariant: the central character and the associated variety. Let us defined them.

If $V$ is an irreducible $\mathfrak{g}$-module, then $End_{U(\mathfrak{g})}(V)=\mathbb{C}$ and thus the elements from $Z(\mathfrak{g})$ act on $V$ as a multiplication by a scalar. Thus we obtain a mapping

$$
\chi: Z(\mathfrak{g}) \rightarrow \mathbb{C}.
$$

It is called {\it the central character}. Actually it depend only on the ideal $Ann(V)$.

To define the associated variety take the graded ring and the graded ideal

$$
grI\subset gr(U(\mathfrak{g}))=S(\mathfrak{g})=\mathbb{C}[\mathfrak{g}]
$$

Then we can take  a Zarisky closed subset in $\mathfrak{g}$ defined by $grI$. It is called {\it the associated variety}. We denote it as $Var(I)$.

Duflo obtained the following classification.

\begin{theorem}Let $\rho$ be half-sum of positive root of the algebra $\mathfrak{g}$, denote as $V(\lambda)$ be an irreducible highest weight module of highest weight $\lambda-\rho$ then
every primitive ideal in $U(\mathfrak{g})$ is of type $Ann(V(\lambda))$ for some weight $\lambda$.
\end{theorem}

It turns out hard to give a criteria when $Ann(V(\lambda))=Ann(V(\mu))$, the solution of this problem is known but very non-trivial.

Now let us proceed to $W_{\chi}$-algebra. Note that by definition $Q_{\chi}$ is a $(U(\mathfrak{g}),W_{\chi})$-bimodule.
 For a $W_{\chi}$-module $V$ we define a $U(\mathfrak{g})$-module $V^+$ as follows

$$
V^{+}=Q_{\chi}\otimes_{W_{\chi}} V.
$$

There exists a map

\begin{equation}
t: Prim W_{\chi}\rightarrow Prim U(\mathfrak{g}), \text{ such that } t(Ann_{W_{\chi}}(V))=Ann_{U(\mathfrak{g})}(V^+)
\end{equation}

\begin{theorem}

This mapping induces a surjection

\begin{equation}
\label{mt}
t: Prim_{fin}W_{\chi}\rightarrow  Prim_{G,e}U(\mathfrak{g}),
\end{equation}
where
\begin{align*}
& Prim_{fin}W_{\chi}=\{I\in Prim W_{\chi}\,\,\, codim I<\infty\},\\
&Prim_{G,e}U(\mathfrak{g})=\{I\in Prim U(\mathfrak{g}), Var(I)=\text{ the closure of }Ge\}
\end{align*}
\end{theorem}

The first set parameterizes  finite dimensional $W_{\chi}$-modules.

Losev obtained some information about fibers of this surjection. Let $C_G(e,f,h)$ be a centralizer in $G$ of elements $e,f,h$. 
  Losev showed that the group $C_G(e,f,h)$ acts in an adjoint way on $U(\mathfrak{g})$ and this action preserves $W_{\chi}$. The identity component $C^0_G(e,f,h)$ acts trivially on $W_{\chi}$ since it's Lie algebra embeds into $W_{\chi}$. Hence, the adjoint action of   $C_G(e,f,h)$  induces an action of $C=C_G(e,f,h)/C^0_G(e,f,h)$ on $W_{\chi}$ and hence on the set $Prim_{fin}W_{\chi}$.

\begin{theorem}
The fibers of \eqref{mt} are orbits of the action of $C$.

\end{theorem}

\subsection{Highest weight theory for $W_{\chi}$}

These results are taken from \cite{bk}. Let us give a definition of a highest weight module. Let us first define a subalgebra $\mathfrak{g}(0)\subset \mathfrak{g}$.

Assume that the toral subalgebra $\mathfrak{t}$ is chosen in such a way that
$$\mathfrak{t}^e=\{t\in\mathfrak{t}:\,\,\,[e,t]=0\}$$

 is maximal toral subalgebra of $\mathfrak{g}^e\cap\mathfrak{g}_0$. Then $\mathfrak{t}^e$ is the orthogonal complement in $\mathfrak{t}$ to $h$.

Here

\begin{align*}
&\mathfrak{g}^e=\{g\in\mathfrak{g}:\,\,\,[e,g]=0\},\,\,\mathfrak{g}_0=\{g\in\mathfrak{g}:\,\,\,[h,g]=0\},\
\end{align*}

For $\alpha\in (\mathfrak{t}^e)^*$ we denote as $\mathfrak{g}_{\alpha}^e$ the $\alpha$-weight space of $\mathfrak{g}^e$. Then

$$
\mathfrak{g}^e=\mathfrak{g}(0)^e\oplus\bigoplus_{\alpha\in \Phi^e}\mathfrak{g}_{\alpha}^e
$$

Here $$\mathfrak{g}(0)^e:=\{g\in\mathfrak{g}^e:\,\,\,\forall h\in\mathfrak{t}^e\,\,\,\, [h,g]=0\}.$$ is a centralizer of $\mathfrak{t}^e$  in $\mathfrak{g}^e$.  For  $\alpha\in (\mathfrak{t}^e)^*$ we put $$\mathfrak{g}_{\alpha}^e:=\{g\in\mathfrak{g}^e:\,\,\,\forall h\in\mathfrak{t}^e\,\,\,\, [h,g]=\alpha(h)g\}.$$ Note that  $\mathfrak{t}^e\subset \mathfrak{g}(0)^e$ but in general  $\mathfrak{t}^e\neq \mathfrak{g}(0)^e$.

Now let us denote the notion of a positive root. Take a Borel subalgebra $\mathfrak{b}\subset\mathfrak{g}$. Take a parabolic subalgebra $\mathfrak{q}=\mathfrak{g}(0)^e+\mathfrak{b}$.

For each $\alpha$ one has either $\mathfrak{g}^e_{\alpha}\subset \mathfrak{q}$ or $\mathfrak{g}^e_{-\alpha}\subset \mathfrak{q}$. In the first case we call the root $\alpha$ positive and in the second case - negative.

Put

\begin{align*}
&\mathfrak{g}_{+}^e=\bigoplus_{\alpha>0}\mathfrak{g}_{\alpha}^e,\\
&\mathfrak{g}_{-}^e=\bigoplus_{\alpha<0}\mathfrak{g}_{\alpha}^e,
\end{align*}
Here $\mathfrak{g}_{\alpha}^e$ is the maximal $\alpha$-weight subspace in $\mathfrak{g}^e$
$\,$
Choose a basis $h_1,...,h_l$ in $\mathfrak{g}^e_0$, $e_1,...,e_m$ in $\mathfrak{g}^e_{+}$,
$f_1,...,f_m$ in $\mathfrak{g}^e_{-}$, introduce notations

\begin{equation}
F_i:=\Theta(f_i),\,\,\, E_i:=\Theta(e_i),\,\,\,\, H_i:=\Theta(h_i)
\end{equation}

We have a PBW basis of $W_{\chi}$ $$F^aH^bE^c.$$ Put \begin{equation*} W^{+}_{\chi}=\text{  a left ideal of $W_{\chi}$ generated by $E_1,...,E_m$}\end{equation*}

Take a $W_{\chi}$-module and for $\lambda\in (\mathfrak{t}^e)^*$ put

\begin{equation}
V_{\lambda}=\{v\in V,\,\,\, tv=\lambda(t)v \text{ for all }t\in\mathfrak{t}^e\}
\end{equation}
 here we use an embedding $\mathfrak{t}^e\subset \mathfrak{g}^e\subset W_{\chi}$.

 \begin{defn}
 We call  $V_{\lambda}$ {\it the weight space}.
 \end{defn}

\begin{defn}
The weight space if {\it maximal} if $W^+_{\chi}V_{\lambda}=0$.
\end{defn}

Now turn to the case of Lie algebras. In this case a highest weight module is a module with a highest weight vector? When in the case of Lie algebras a module with the maximal weight space has a highest weight vector? The answer is the following: when this maximal weight space is one-dimensional and it generates the whole module. We can reformulate the condition that this maximal weight space is one-dimensional as follows: it is an irreducible representation of the Cartan subalgebra.
We use this reformulation for the definition of the  maximal weight modules of $W_{\chi}$.

Note that $e\in\mathfrak{g}(0)$, thus we can consider the $W$-algebra  $W^{\mathfrak{g}(0)}_{\chi}$. From the second definition it immediately follows that  $W^{\mathfrak{g}(0)}_{\chi}\subset W_{\chi}$. In the case of finite $W$-algebras the role of Cartan subalgebra is played by $W^{\mathfrak{g}(0)}_{\chi}$.

\begin{defn}
A $W_{\chi}$-module is  {\it a module of maximal weight} if  it is generated by maximal weight space, which is finite dimensional and irreducible as a $W^{\mathfrak{g}(0)}_{\chi}$ - module.
\end{defn}

\begin{defn}
For a finite dimensional and irreducible $W^{\mathfrak{g}(0)}_{\chi}$ - module $V$ we  define a Verma module

$$
M= W_{\chi}/(W_{\chi} \otimes_{W^{    \mathfrak{g}(0)}_{\chi}  } V)
$$
\end{defn}
This module has a unique maximal proper submodule $R$ and $M/R$ is an irreducible module.

Every irreducible module can be obtained in this way.

\subsection{The category of representations}

The following results about the category of representations of $W_{\chi}$ are known

\begin{theorem}
The number of isomorphism classes of  irreducible finite-dimensional representations with a given central character is finite
\end{theorem}

\begin{theorem}
Every Verma module $M$ has a finite composition series
\end{theorem}

Now take a category $\mathcal{O}$ - the category of finitely generated $W_{\chi}$-modules that are semisimple over $\mathfrak{t}^e$ with finite-dimensional $\mathfrak{t}^e$-weight-spaces and such that the set $\{\lambda\in (\mathfrak{t}^e:\,\,\, V_{\lambda}\neq 0)\}$ is contained in $\{\nu\in (\mathfrak{t}^e)^*:\,\,\,\, \nu\leq \mu\}$

\begin{theorem}
Every object in $\mathcal{O}$ has a composition series and the category $\mathcal{O}$ decomposes into a direct sum of categories $\mathcal{O}_{\psi}$, where $\mathcal{O}_{\psi}$
is a category generated by irreducible modules with a central character $\psi$
\end{theorem}

\section{A relation to shifted yangians}

\subsection{The shifted Yangians and finite $W$-algebras}

In \cite{bk} (see also \cite{bk1}) the following result was obtained. A yangian $Y(gl_n)$ is an algebra generated by elements $1$, $t_{i,j}^r$, $r>0$, $i,j=1,...,N$ subject to relations

\begin{equation}
[t_{i,j}^{r+1},t_{k,l}^s]-[t_{i,j}^r,t_{k,l}^{s+1}]=t_{k,j}^rt_{i,l}^s-t_{k,j}^st_{i,l}^r.
\end{equation}

Take a series

$$
t_{i,j}(u)=\sum_{r\geq 0}t_{i,j}^ru^{-r}
$$

A matrix $T(u)=(t_{i,j}(u))$ can be considered as an element of the space

$$
T(u)=\sum E_{i,j}\otimes t_{i,j}(u)\in gl_{n}\otimes Y(gl_n)
$$

Define an element

$$
R(z)=1-\frac{\sum E_{i,j}\otimes E_{j,i}}{z}
$$

Then in the space $  gl_{n}\otimes  gl_{n}\otimes Y(gl_n)$ the defining commutation relations of the Yangian can be written as follows

$$
R_{1,2}(u-v)T_{1,3}(u)T_{2,3}(v)=T_{2,3}(v)T_{1,2}(u)R_{1,2}(u-v)
$$

Consider  a Gauss decomposition

$$
T(u)=F(u)D(u)E(u),
$$

where

\begin{align*}
&D(u)=diag(D_1(u),...,D_n(u),\,\,\, E(u)=\begin{pmatrix}   1 & E_{1,2}(u)&...& E_{1,n}(u)\\
0& 1&...& E_{2,n}(u)\\
...\\
0&0&...& 1 \end{pmatrix}\\
& F(u)=\begin{pmatrix}   1 & 0&...& 0\\
E_{2,1}(u)& 1&...& 0\\
...\\
E_{n,1}(u)&E_{n,2}(u)&...& 1 \end{pmatrix}
\end{align*}
Put

$$
\tilde{D}_i(u):=-D_i(u)^{-1}
$$

One can decompose the matrix elements of these matrices

\begin{align*}
&E_{i,j}(u)=\sum_{r\geq 1}E_{i,j}^ru^{-r},\,\,\,\,\,\,F_{i,j}(u)=\sum_{r\geq 1}F_{i,j}^ru^{-r}\\
&D_{i}(u)=\sum_{r\geq 1}D_{i}^ru^{-r},\,\,\,\,\,\,\tilde{D}(u)=\sum_{r\geq 1}\tilde{D}^ru^{-r}\\
\end{align*}

 and obtain another set of generators of $Y(gl_N)$. The commutation relations for them are the following:
\begin{align*}
&[D_i^r,D_j^r]=0,\,\,\,\ [E_i^,F_j^s]=\delta_{i,j}\sum_{t=0}^{r+s-1}\tilde{D}_i^tD_{i+1}^{r+s-1-t}\\
&[D_i^r,F_j^s]=(\delta_{i,j}-\delta_{i,j+1})\sum_{t=0}^{r-1}F_j^{r+s-1-t}D_i^t\\
&[D_i^r,E_j^s]=(\delta_{i,j}-\delta_{i,j+1})\sum_{t=0}^{r-1}D_i^tE_j^{r+s-1-t}\\
&[E_i^r,E_i^s]=\sum_{t=1}^{s-1} E_i^tE_i^{r+s-1-t}-\sum_{t=1}^{s-1} E_i^{r+s-1-t}E_i^t\\
&[F_i^r,F_i^s]=\sum_{t=1}^{s-1} F_i^tF_i^{r+s-1-t}-\sum_{t=1}^{s-1} F_i^{r+s-1-t}F_i^t
\end{align*}

\begin{align*}
&[E_i^r,E_{i+1}^{s+1}]-[E_i^{r+!},E_{i+1}^s]=-E_{i}^rE_{i+1}^s\\
&[F_i^r,F_{i+1}^{s+1}]-[F_i^{r+!},F_{i+1}^s]=-F_{i}^rF_{i+1}^s\\
&[E_{i}^r,E_j^s]=[F_{i}^r,F_j^s]=0\,\,\,|i-j|>1\\
&[E_i^r,[E_i^s,E_j^t]]+[E_i^s,[E_i^r,E_j^t]]=0\,\,\,|i-j|>1\\
&[F_i^r,[F_i^s,F_j^t]]+[F_i^s,[F_i^r,F_j^t]]=0\,\,\,|i-j|>1\\
\end{align*}

Take a matrix $S=(s_{i,j})$, such that

\begin{equation}
s_{i,j}+s_{j,k}=s_{i,k}\text{ if }|i-j|+|j-k|=|i-j|.
\end{equation}

For a Young  diagram $\lambda=(\lambda_1,...,\lambda_n)$ we can define  a shift matrix $S_{\lambda}$ by formula

$$
s_{i,j}=\lambda_{n+1-j}-\lambda_{n+1-i}\text{ for $i<j$ and  $0$ if $i\geq j$ }
$$

\begin{defn}

A {\it shifted Yangian} is a subalgebra in $Y(\mathfrak{gl})_n$ generated by

\begin{align*}
& D_i^r,\,\,\, r>0,\,\, 1\leq i\leq n\\
& E_{i,i+1}^r\,\,\,\, r>s_{i,i+1}\,\,\, 1\leq i<n\\
& F_{i+1,i}^r\,\,\,\, r>s_{i+1,i}\,\,\, 1\leq i<n\\
\end{align*}

\end{defn}

To $\lambda$ one assigns a nilpotent element $e_{\lambda}\in \mathfrak{gl}_N$. Denote the finite $W$ algebra corresponding to this nilpotent element as  $W_{e_{\lambda}}$.

\begin{theorem}
$Y(gl_N,S_{\lambda})=W_{e_{\lambda}}$
\end{theorem}

Similar results are know for the series $B$, $D$, $C$ (see \cite{br1}, \cite{br2})

\subsection{Representations of shifted yangians}

The theory of finite-dimensional representations of shifted yangiangs is well-developed.

\begin{prop}
Every irreducible  finite-dimensional module has a vector $v_+$ annihilated by $E_{i,i+1}^r$ on which $D_i^r$ act diagonally. The isomorphism type of a module is defined by $D_i^r$-eigenvalues,  for each $i$ these eigenvalues are zero for $r>\lambda_r$.
\end{prop}

Let us write

\begin{align*}
& u^{p_1}D_1(u)v=P_1(u)v\\
& (u-1)^{p_2}D_2(u-1)v=P_2(u)v\\
&...\\
&(u-n+1)^{p_n}D_n(u-n+1)v=P_n(u)v,
\end{align*}

where $P_1,...,P_n$ are monic\footnote{i.e. the coefficient at the highest power of  $u$ equals to $1$} polynomials of degrees $\lambda_1,...,\lambda_n$.

The corresponding module we denote $V(P_1(u),...,P_n(u))$.

Put
$$
P_i(u)=(u+a_{i,1})...(u+a_{i,\lambda_i}),\,\,\, a_{i,1}\leq...\leq a_{i,\lambda_i}.
$$

The considered shifted yangian is defined by shifted Young diagram $\lambda=(\lambda_1\geq...\geq \lambda_n)$. We can represent it as usual as a collection of boxes: $\lambda_1$ boxes in the upper row,  $\lambda_2$ boxes in the next row and so on.

Put elements  $a_{i,j}$ the diagram and obtain a Young tableau.

\begin{theorem}
This module is irreducible if and only if  tableau is standard (the entries increase from bottom to columns)
\end{theorem}

\section{A relation to the center of the universal enveloping algebra}


\subsection{Sugawara elemnts}
It turns out that in general the center of $U'(\widehat{\mathfrak{g}})_k$ is trivial.

\begin{theorem}[see  \cite{edfr}]
If $k+g\neq 0$ then the center of $U'(\widehat{\mathfrak{g}})_k$ is generated by $1$.
\end{theorem}

\begin{defn}[see \cite{Ha}]
We say that an element $P\in U'(\widehat{\mathfrak{g}})_k$ is a {\it Sugawara element} if

$$
[P,U'(\widehat{\mathfrak{g}})_k]=(k+g)U'(\widehat{\mathfrak{g}})_k
$$
\end{defn}

Let us give an example of a Sugawara element. Take a Casimir of the second order in $U'(\mathfrak{g})_k$:

\begin{equation}
T=\sum_{\alpha} I_{\alpha} I^{\alpha}
\end{equation}

Let us write field

\begin{equation}
T(z)=\sum_{\alpha} (I_{\alpha} I^{\alpha})(z)=\sum z^{-n-2}\bar{L}_n
\end{equation}

Note that we put no constant in front of the sum in this formula. Thus these elements are defined for all $k$ and for $k\neq -g$ we have $\bar{L}_n=(k+g)L_n$.

\begin{prop}
The elements $\bar{L}_n$ are Sugawara elements
 \end{prop}

It turns out the the generalization of this result to central elements of higher orders is not direct. If we take an central element in $U'(\mathfrak{g})$ of higher order

$$
T_n=d_{\alpha_1,...,\alpha_n}I_{\alpha_1}...I_{\alpha_n}
$$

then the modes of the field
$$
T_n(z)=d_{\alpha_1,...,\alpha_n}(I_{\alpha_1}(...I_{\alpha_n})...)(z)=\sum_{k}z^{-n-k}T_n^k
$$

are not in general the Sugawara elements.

Nevertheless the following result is know.

\begin{theorem}[see  \cite{edfr}]
There exist Sugawara elements $\bar{T}_n^k$ such that

$$
\bar{T}_n^k=T_n^k\,\,\,mod F_{k-1}U'(\widehat{\mathfrak{g}})
$$

 The center is a polynomial algebra generated by elements $\bar{T}_n^k$.

\end{theorem}

Note that for all values of $k$  the Sugawara elements form an associative algebra. For $k=-g$ it is of course commutative. Now put $h=k+g$, then we obtain a family of  associative algebras depending on $h$, such that for $h=0$ the corresponding algebra is commutative.  Look at this family from the point of view of Section \ref{defq}. We have a deformation of a commutative algebra for $h=0$ and this deformation must be described by a Poisson bracket. Thus we come to the following result

\begin{lem}
The algebra of   Sugawara elements for $k=-g$ has a natural  structure of a Poisson algebra.
\end{lem}

\subsection{Talalaev's construction}
For a long time an explicit construction of  the elements $\bar{T}_n^k$ was unknown. But recently it was given  for  $\mathfrak{g}=\mathfrak{gl}_N$ by Talalaev.

Let us remind a construction of generators of the center of $U(\mathfrak{gl}_N)$ using determinants. Let $\tau$ be a formal variable, consider the matrix

\begin{equation}
\Omega=\begin{pmatrix} E_{1,1}+\tau+1 & E_{1,2} &...& E^1_{1,N}\\ E^1_{2,1} & E^1_{2,2}+\tau+2 &...& E^1_{2,N}\\...\\ E^1_{N,1}& E^1_{N,2}&...& E^1_{N,N}+\tau+N\end{pmatrix}.
\end{equation}

Take it's row determinant

$$
rdet\Omega=\sum_{\sigma\in S_N}(-1)^{\sigma  }\Omega_{1,\sigma(1)}...\Omega_{N,\sigma(N)},
$$

And consider it's decomposition

$$
rdet\Omega(u)=\tau^N+c_1\tau^{N-1}+...+c_N
$$

\begin{theorem}
The elements $ c_i$ freely generate the center of $U(\mathfrak{gl}_N)$
\end{theorem}

Now consider the Talaev's construction of Sugawara elements. Consider $\widehat{\mathfrak{gl}_N}$, it is spanned by $C$ and $E_{i,j}^n$.
 Introduce a notation
$$
\tau(E_{i,j}^n):=E_{i,j}^{n+1}
.$$

Take a matrix

\begin{equation}
\Omega=\begin{pmatrix} E^{-1}_{1,1}+\tau & E^{-1}_{1,2} &...& E^{-1}_{1,N}\\ E^{-1}_{2,1} & E^{-1}_{2,2}+\tau &...& E^{-1}_{2,N}\\...\\ E^{-1}_{N,1}& E^{-1}_{N,2}&...& E^{-1}_{N,N}+\tau\end{pmatrix}.
\end{equation}

Take it's determinant and it's decomposition

\begin{equation}
rdet\Omega=\tau^N+c_1\tau^{N-1}+...+c_N
\end{equation}

\begin{theorem}[see \cite{tal} or \cite{mc}]
The elements  $c_i^k:=\tau^k c_i$, $k\in\mathbb{Z}$ are Sugawara elements and they generate the algebra of Sugawara elements.
\end{theorem}

In another way we can describe this construction as follows.   Consider a matrix composed of currents

\begin{equation}
\Omega(u)=\begin{pmatrix} E_{1,1}(z)+\frac{d}{dz} & E_{1,2}(z) &...& E_{1,N}(z)\\ E_{2,1}(z) & E^1_{2,2}(z)+\frac{d}{dz} &...& E_{2,N}(z)\\...\\ E_{N,1}(z)& E_{N,2}(z)&...& E_{N,N}(z)+\frac{d}{dz}\end{pmatrix}.
\end{equation}

Then take a determinant using a normal ordered product

\begin{equation}
rdet\Omega=\sum_{\sigma\in S_N}(-1)^{\sigma  }(\Omega_{1,\sigma(1)}(...\Omega_{N,\sigma(N)}))(u)=\frac{d^N}{dz^N}+c_1\frac{d^{N-1}}{dz^{N-1}}+...+c_N
\end{equation}

Finally consider decompositions

\begin{equation}
c_k(z)=\sum_k c_i^k z^{-i-k}.
\end{equation}

Then we obtain just the elements  $c_i^k$ mentioned in the theorem above.

Analogous construction were later given for the series $B,C,D$ and for $G_2$  (see \cite{m1}, \cite{m2}, \cite{m3}).

Remind that the algebra $W_N$ was defined using the central elements in the universal enveloping algebra   $U(\mathfrak{gl}_N)$. The Sugawara elements are closely related to this center. Thus the following result is not surprising.

\begin{theorem}
The Poisson algebra of Sugawara elements for $\mathfrak{gl}_N$ and $k=-g$ is isomorphic to the classical $W_N$ algebra.
\end{theorem}

There are analogs of this theorem for other series of algebras. The Poisson algebra of Sugawara elements for  the algebra  $\mathfrak{g}$ is isomorphic to the classical $W$-algebra associated to a principle nilpotent element but for the Langlands dual Lie algebra, i.e. Lie algebra with a transpose  Cartan matrix.

 \subsection{Explicit generators of the $W_N$}
 \label{wnm}

Actually the technique presented above allows to obtain generators of the algebra $W_N$ (see \cite{m4}).

Put

\begin{equation}
\Omega=\begin{pmatrix} E^{-1}_{1,1}+N\tau & E^{-1}_{1,2} & E^{-1}_{1,3}  &...& E^{-1}_{1,N-1}& E^{-1}_{1,N}\\ -1 & E^{-1}_{2,2}+N\tau & E^{-1}_{2,3}  &...& E^{-1}_{2,N-1}& E^{-1}_{2,N}\\ 0 & -1 &E^{-1}_{3,3}+N\tau &...& E^{-1}_{3,N-1}& E^{-1}_{3,N}\\...\\ 0& 0&...&-1& E^{-1}_{N,N}+N\tau\end{pmatrix}.
\end{equation}

Then we take the determinant and it's decomposition

\begin{equation}
rdet\Omega=\tau^N+W_1\tau^{N-1}+...+W_N
\end{equation}

\begin{theorem}[\cite{m4}]
The elements $W_i^k:=\tau^k W_i$, $k\in\mathbb{Z}$ belong to $W_N$ and generate it as an associative algebra.
\end{theorem}

\subsection{The case of a finite $W$-algebra}

\label{cent}

Let us formulate another result about the structure of a finite $W$-algebra

\begin{theorem}[\cite{ko}]
\label{tt}
If $e$ is a regular nilpotent element then $W_{\chi}$ is isomorphic to the center of $U(\mathfrak{g})$.
\end{theorem}

\subsubsection{Examples}

Let us describe this isomorphism explicitly in some examples. Below we give an explicit description of the $W$-algebras associated with a principle nilpotent elements. We give generators corresponding to generators of the center of $U(\mathfrak{g})$.

Take $\mathfrak{g}=\mathfrak{gl}_2$, then $B(a,b)=tr(ab)$, take the $\mathfrak{sl}_2$-triple

\begin{align*}
&e=
\begin{pmatrix}
0 &1\\
0 & 0
\end{pmatrix}
f=
\begin{pmatrix}
0 &0\\
1 & 0
\end{pmatrix}
h=
\begin{pmatrix}
1 &0\\
0 & -1
\end{pmatrix}
\end{align*}

The Dynkin grading is given by the matrix

\begin{align*}
&\begin{pmatrix}
 0& 2\\ -2 & 0
\end{pmatrix}
\end{align*}

Take $z=\begin{pmatrix}   1& 0\\ 0&1\end{pmatrix}$, then

\begin{align*}
&\mathfrak{g}^e=<z,e>,\\
&\mathfrak{m}=\mathfrak{g}_{-2}=<f>,\,\,\,\,\chi(f)=B(f,e)=1,\\
&\mathfrak{v}=\mathfrak{g}_0\oplus\mathfrak{g}_2
\end{align*}

The algebra $W_{\chi}$ is freely generated by elements $\bar{x}$, $\overline{e+\frac{1}{4}h^2+\frac{1}{2}}=\overline{\Omega}$, where $\Omega=ef+fe+\frac{1}{2}h^2$ is the second Casimir element.  Here as in Section \ref{sedef} we denote as $\bar{x}$ the projection $U(\mathfrak{g})\rightarrow U(\mathfrak{g})/I_{\chi}$.

The center of $U(\mathfrak{gl}_2)$ is just   freely  generated by $z$ and $\Omega$. This is an illustration of the isomorphism form Theorem \ref{tt}.



Take now $\mathfrak{g}=\mathfrak{gl}_N$, $B(a,b)=tr(ab)$,

\begin{align*}
&e=\begin{pmatrix}   0 & 1  & 0 & ... &0\\  0&0  &1 &...&0\\...\\ 0&0&0&...&1\\
0&0&0&...&0\end{pmatrix},\,\,\,
f=\begin{pmatrix}   0 & 0  & 0 & ... & 0&0\\  0&n-1  &0 &...&0&0\\ 0&0  &2(n-2) &..&0.&0\\...\\ \
0&0&0&...& n-1&0\end{pmatrix}\\
& h=diag(n-1,n-3,...,3-n,1-n).
\end{align*}

The Dynkin grading is given by the matrix

\begin{align*}
\begin{pmatrix}
0&2&4& 6&...&2n-2\\
-2&0&2&4&...&2n-4\\
-6&-4&-2&0&...&2n-6\\
...\\
2-2n&&&&...&0
\end{pmatrix}
\end{align*}

As before $z=diag(1,....1)$, then

\begin{align*}
&\mathfrak{g}^e=<z,e,e^2,...,e^{n-1}>\\
&\mathfrak{m}=\oplus_{j\geq 2}\mathfrak{g}_{2-2j}, \,\,\,\, \chi(E_{i+1,j})=1,\,\,\,\chi(E_{i+k,j})=0\text{ for }k\geq 2
\end{align*}

Take Casimir elements

$$
\Omega_k=E_{i_1,i_2} E_{i_2,i_3}...E_{i_{k}i_1}
$$
 Then $W_{\chi}$ is a polynomial algebra generated by

$$
\bar{z},\,\,\,\overline{\Omega_k}
$$

But the center of $U(\mathfrak{gl}_N)$ is just   freely  generated by $z$ and $\Omega_k$. So we again have  an illustration of the isomorphism form Theorem \ref{tt}.


Unfortunately the descriptions of $W_{\chi}$ obtained above are not explicit sice we describe the generators as projections in the factor algebra $U(\mathfrak{g})/I_{\chi}$ as in the first definition of $W_{\chi}$. But due to the second definition of finite $W$-algebra  we have an embedding $W_{\chi}=U(\mathfrak{v})^{\mathfrak{m}}\subset U(\mathfrak{g})$.  Let us give explicit  formulas for the generators of $W_{\chi}$  as elements of  $U(\mathfrak{g})$. This description is similar to description of generators of infinite $W$-algebra obtain in the previous Section.

Consider the case
 $\mathfrak{g}=\mathfrak{gl}_N$, $e=E_{1,2}+E_{2,3}+...+E_{N-1,N}$, let

\begin{equation}
\Omega(u)=\begin{pmatrix}
  E_{1,1}+u-1 & E_{1,2} & E_{1,2} &...& E_{1,n}\\
  1& E_{2,2}+u-2  & E_{2,3}&...&E_{2,n}\\
 0&1 & E_{3,3}+u-3&...& E_{3,n}\\
...\\
0&0&0&...&E_{n,n}+u-n
   \end{pmatrix}
\end{equation}

Consider the row determinant

$$
rdet\Omega(u)=\sum_{\sigma\in S_N}(-1)^{\sigma}\Omega_{1,\sigma(1)}...\Omega_{N,\sigma(N)}.
$$

Take a decomposition $$rdet\Omega(u)=u^n+\sum_{i=1}^n w_i u^{n-i}$$

\begin{theorem}\cite{bk}
The elements $w_i$ belong to $W_{\chi}\subset U(\mathfrak{g})$, commute and freely generate $W_{\chi}$
\end{theorem}

\section{Noncommutative pfaffians and Capelli elements}

\subsection{The definition}

Take an algebra $\mathfrak{o}_N$ and consider it's generators $F_{ij}=E_{ij}-E_{ji}$, $i,j=1,...,n$, $i<j$. They satisfy the relations

\begin{equation}
[F_{ij},F_{kl}]=\delta_{kj}F_{il}-\delta_{il}F_{kj}-\delta_{ik}F_{jl}+\delta_{jl}F_{ki}.
\end{equation}

Let  $\Phi=(\Phi_{ij})$, $i,j=1,...,2k$  be a skew-symmetric
 $2k\times 2k$-matrix whose elements belong to some ring.

A noncommutative pfaffian of  $\Phi$ is defined by formula

\begin{equation}
Pf \Phi=\frac{1}{k!2^k}\sum_{\sigma\in
S_{2k}}(-1)^{\sigma}\Phi_{\sigma(1)\sigma(2)}...\Phi_{\sigma(2k-1)\sigma(2k)},
\end{equation}

Take a matrix $F=(F_{ij})$, $i,j=1,...,N$.  Since
$F_{ij}=-F_{ji}$ then $F$ is skew symmetric.

Take a set of indices
$I\subset\{1,...,N\}$ and consider a matrix

$$
F_I=(F_{i,j})_{i,j\in I}
$$

\begin{theorem} [\cite{Molev}]
Let

\begin{equation}
C_k=\sum_{|I|=k,I\subset\{1,...,N\}}(PfF_I)^2,\,\,\,\,k=2,4,...,2[\frac{N}{2}].\end{equation}
 Then $C_k$  belong to the center of $U(\mathfrak{o}_N)$. In the case of odd $N$ they generate the center. In the case of odd
 $N$ they generate the center if we add $PfF$.
\end{theorem}

In \cite{AG} we found commutation relations between $PfF_{I}$ and the generators. Introduce a notation $F_{ij}I$.Let  $I=\{i_1,...,i_k\}$, $i_r\in\{1,...,N\}$
be a set of indices. Identify
$i_r$ with   the vector $e_{i_r}$, and the set $I=\{i_1,...,i_k\}$ we identify with
 $e_{i_1}\otimes ...\otimes e_{i_k}\in V^{\otimes k}$. Then $F_{ij}I$  is a result of the action of
  $F_{ij}$ onto $I$.

For $\alpha,\beta\in\mathbb{C}$ put  $$PfF_{\alpha I+\beta
J}:=\alpha PfF_I+\beta Pf F_J.$$ Then for every
$g\in\mathfrak{o}_N$ we have defined  $PfF_{gI}$.

\begin{prop}
\label{pffco} $[F_{ij},PfF_I]=PfF_{F_{ij}I}$
\end{prop}

Also we found relations between two pfaffians

\begin{prop}
$$
[PfF_I,PfF_{J}]=\sum_{I=I'\cap I''}(-1)^{(I'I'')}PfF_{I'}PfF_{PfF_{I''}J}
$$
\end{prop}

These relations allow to prove the first part of the Molev's theorem.

\begin{conjecture}
What kind of algebraic object form pfaffians? Is it isomorphic to a finite $W$ algebra?
\end{conjecture}

\subsection{Fields corresponding to Capelli elements}

Define a field $PfF_I(z)$, $I=\{i_1,...,i_{2k}\}$  as follows (see \cite{AG2}):

\begin{defn}

\begin{align}
&Pf F_{I}(z)=\frac{1}{k!2^k}\sum_{\sigma\in
S_{2k}}(-1)^{\sigma}(F_{\sigma(i_1)\sigma(i_2)}(F_{\sigma(i_3)\sigma(i_4)}(F_{\sigma(i_5)\sigma(i_6)}...\\&...F_{\sigma(i_{2k-1})\sigma(i_{2k})}))...)(z).
\end{align}

\end{defn}

Introduce a field by analogy with a Capelli element

\begin{equation}
C_n(z)=\sum_{|I|=n}(PfF_IPfF_I)(z).
\end{equation}

In this case we have the following analogue of the formula  \eqref{pffco}

\begin{align}
\begin{split}
\frac{(c_{n}k+d_{n})PfF_{I\setminus J}(w)}{(z-w)^2}+\frac{PFF_{F_J}I(w)}{(z-w)},
\end{split}
\end{align}

where
\begin{align}
\begin{split}
&c_{n}k+d_{n}=k\frac{2}{n}+\frac{4(n-1)}{n}+(c_{n-1}k+d_{n-1})\frac{2}{n}\frac{(n-1)(n-2)}{2}=\\
&=\frac{1}{2}(2+c_{n-1}(n-1)(n-2))k+\frac{1}{n}((4n-4)+c_n(n-1)(n-2).
\end{split}
\end{align}

\begin{conjecture}
Are the modes of elements $C_n(u)$ Sugawara elements? If no then how should we change them to obtain Sugawara elements?
\end{conjecture}

Actually the following result takes place

\begin{prop}
\begin{align*}
&[C_k^n,I_{\alpha}^s]=(k+g)\text{ for } k=-n-1,-n-2,...,\\
&[C_k^n,I_{\alpha}^s]=0\text{ for } k=0,1,...
\end{align*}
\end{prop}

\end{document}